\documentclass[12pt]{article}
\usepackage{amsfonts}
\usepackage{mathrsfs}
\usepackage{txfonts}
\newtheorem{prop}{Proposition}
\newtheorem{cor}{Corollary}
\newtheorem{lem}{Lemma}
\newtheorem{thm}{Theorem}
\newtheorem{defi}{Definition}

\newcommand {\beq}{\begin{equation}}
\newcommand {\eeq}{\end{equation}}
\newcommand {\beqa}{\begin{eqnarray}}
\newcommand {\eeqa}{\end{eqnarray}}         %Equation version
\newcommand {\beqas}{\begin{eqnarray*}}
\newcommand {\eeqas}{\end{eqnarray*}}
\newcommand {\bea}{\begin{array}}
\newcommand {\eea}{\end{array}}
\newcommand {\bds}{\begin{displaymath}}
\newcommand {\eds}{\end{displaymath}}

\newcommand {\nn}{\nonumber}
\newcommand{\no}{\noindent}
\newcommand {\bebb}{\begin{thebibliography}}
\newcommand {\eebb}{\end{thebibliography}}      %Reference version
\newcommand {\bbit}{\bibitem}
\def\dl{\delta}

\def\journal#1&#2(#3){\unskip, \sl #1\ \bf #2 \rm(19#3) }
\def\andjournal#1&#2(#3){\sl #1~\bf #2 \rm (19#3) }

\begin{document}
%\begin{titlepage}

\begin{flushright}
\end{flushright}

\baselineskip =18pt

\vskip 1cm

\begin{center}
%\title
{\Large\bf The elliptic quantum algebra $U_{q,p}(\widehat{sl_N})$
and its vertex operators}

\vspace{1cm}

%\author

{\normalsize\bf Wen-Jing Chang $^{a, b}$ and Xiang-Mao Ding
$^a$\footnote{corresponding author; E-mail: xmding@amss.ac.cn}}
\vspace{1cm}

{\em $^1$ Institute of Applied Mathematics, Academy of Mathematics
and Systems Science; \\Chinese Academy of Sciences, P.O.Box 2734,
Beijing 100190, People's Republic of China}
\\
{\em $^2$ Graduate School of Chinese Academy of Sciences, \\
Beijing 100049, People's Republic of China}

\end{center}
\date{}

%\maketitle
\vspace{2cm}

\begin{abstract}

We construct a realization of the elliptic quantum algebra
$U_{q,p}(\widehat{sl_N})$ for any given level $k$ in terms of free
boson fields and their twisted partners. It can be considered as the
elliptic deformation of the Wakimoto realization of the quantum
affine algebra $U_{q}(\widehat{sl_N})$. We also construct a family
of screening currents, which commute with the currents of
$U_{q,p}(\widehat{sl_N})$ up to total q-differences. And we give
explicit twisted expressions for the type $I$ and the type $II$
vertex operators of $U_{q,p}(\widehat{sl_N})$ by twisting the known
results of the type $I$ vertex operators of the quantum affine
algebra $U_{q}(\widehat{sl_N})$ and the new results of the type $II$
vertex operators of $U_{q}(\widehat{sl_N})$ we obtained in this
paper.

\end{abstract}
\vspace{1cm}

%%%%%PACS: 11.25Hf; 11.30.Rd; 03.65Fd; 02.20.Hj.

\vspace{0.5cm}

%\end{titlepage}

%\setcounter{section}{1}
%\setcounter{equation}{0}

\section{Introduction}
\vskip 0.5cm

Infinite-dimensional symmetries, such as the Virasoro algebra
($W$-algebra in more general) and affine Lie algebras play central
roles in the two-dimensional conformal field theories (2D CFT)
\cite{CFT}. For the non-conformal (off-critical) integrable
theories, their roles are taken over by the so called quantum
algebras. From the algebraic point of view, there are three kinds of
quantum algebras, according to different exchange properties, which
are nominated as rational, trigonometric and elliptic quantum
algebras respectively. The quantum algebras of the former two kinds
could be regarded as certain degenerate cases of the latter one. For
example, the quantum affine algebras (trigonometric), which are also
known as the quantum group \cite{Drinfeld86,Jimbo}, and the Yangian
double \cite{Smir} with central (rational) can be obtained as a
certain limited case of the elliptic quantum algebras. Various
versions of elliptic quantum algebras, also called as elliptic
quantum groups \cite{Felder,Fronsdal,EF} have been introduced to
understand elliptic face models of statistical mechanics, and in
their semiclassical limit, CFT of Wess-Zumino-Witten (WZW) models on
tori. Their roles are similar to the Kac-Moody algebras in WZW
models. From the Hopf algebra point of view, the elliptic quantum
groups are nothing but quantum affine algebras equipped with a
co-product different from the original one by a certain kind of
twisting, so they can be viewed as quasi-Hopf algebras in the sense
of Drinfeld \cite{Drinfeld}. They have two types which correspond to
different types of integrable models: the vertex type
$A_{q,p}(\widehat{sl_N})$ and the face type
$B_{q,\lambda}(\mathcal{G})$, where $\mathcal{G}$ is a Kac-Moody
algebra associated with a symmetrizable generalized Cartan matrix
\cite{Kac}. The former is closely related to vertex models, for
example, the XYZ model, or equivalently, the eight vertex model in
the principal regime \cite{Baxter}; while some face models, such as
the Andrew-Baxter-Forrester (ABF) models \cite{ABF} which are
`solid-on-solid' (SOS) face models, possess symmetries corresponding
to the face type elliptic algebras $B_{q,\lambda}(\mathcal{G})$.

In mathematics, it is natural to study these algebraic objects'
structures and their representations. In physical applications,
their representations are also required. The standard scheme to
study integrable models in field theories or statistical mechanics
is solving the following basic problems: to diagonalize the given
Hamiltonian and then to compute the correlation functions. Usually,
it is quite difficult to solve such problems directly. It has been
indicated that the algebraic analysis method is an extremely
powerful tool in studying solvable lattice models, especially in
deriving the correlation functions. This method is based on the
infinite dimensional quantum group symmetry possessed by a solvable
lattice model and the representation theory of such symmetry. This
algebraic method could be viewed as the quantum version of the
powerful Inverse Scattering Method \cite{FaddTak}. See \cite{JimMiw}
for a review on it. As a result, if one expects to perform algebraic
analysis over the above two types of elliptic lattice models, he
should first study the corresponding elliptic quantum groups and
their representations.

It is of special interest for the algebra of the intertwining
operators in the WZW model. It was derived by Knizhnik and
Zamolodchikov that the matrix coefficients of the intertwining
operators for the WZW model satisfy certain holonomic differential
equations, i.e., the Knizhnik-Zamolodchikov(KZ) equation \cite{KZ}.
In \cite{FR}, for quantum affine algebra, the authors defined
q-deformed vertex operators as certain intertwining operators and
showed that they satisfied some holonomic difference equations
called the quantum KZ(qKZ) equations. So it is also expected that
the representations of the elliptic quantum algebras are helpful in
constructing the elliptic type solutions of quantum
Knizhnik-Zamolodchikov-Bernard(qKZB) equation, which is a higher
genus extension of the qKZ equation \cite{Bern}.

At classical level, there are various models of representations for
the current algebras and each of them is of significance in certain
applications. Here, we just mention two of them: the Wakimoto
construction (free field realization) \cite{Wak,FF,Yu}, and the
parafermion realization \cite{Nemesch,GepQiu,Gep}. Recently, the
explicit description of free field realizations of current algebras
were given in \cite{Feh1,Feh2,YangZh}. In \cite{DFJMN, JMMN}, the
XXZ model in the anti-ferromagnetic regime was solved by applying
the level one representation theory of the quantum affine algebra
$U_q(\widehat{sl_2})$. In studying a higher spin extension of the
XXZ model, the realizations of $U_q(\widehat{sl_2})$ at level $k>1$
are required, and they were constructed by several authors, such as
the Wakimoto realization in \cite{Shiraishi} and the parafermion
realizations in \cite{Matsuo, DingW}. Furthermore, in \cite{AOS},
free field realization of $U_q(\widehat{sl_N})$ with arbitrary level
$k\geq 1$ was given, and it plays a central role in understanding
the higher rank extension of the XXZ model. The Wakimoto
construction is also a powerful way to study the integrable massive
field theories \cite{Luky}. In practice, free field realization,
which is an infinite dimensional extension of the Heisenberg
algebra, is quite an effective and useful approach to study
complicated algebraic structures and their representations. The
level $k$ free field representation of Yangian double
$DY_{\hbar}(sl_2)$ and applications in physical problems were
discussed in \cite{IK, Konno2}. The level one free field realization
of the Yangian double with central $DY_{\hbar}(sl_N)$ was
constructed in \cite{Iohara}, while the level $k$ representation of
$DY_{\hbar}(sl_N)$ and $DY_{\hbar}(gl_N)$ were given in \cite{DHHZ}.
It should also be remarked that the Yangian double with central
$DY_{\hbar}(\widehat{sl_2})$ is the symmetry possessed by the
Sine-Gordon model, which is the field theory limit of the restricted
SOS (RSOS) model \cite{DHZ,HY}.

It is first noticed by Lukyanov and Pugai \cite{Lukpugai} that a
symmetry of the RSOS model is generated by the q-deformation of the
Virasoro algebra (q-Virasoro algebra). The free field realizations
of screening currents and vertex operators enable them to analyze
the structure of the highest weight representation of the q-Virasoro
algebra. And the screening currents they constructed satisfy an
elliptic deformation of $U_q(\widehat{sl_2})$ at level one, which is
called the elliptic algebra $U_{q,p}(\widehat{sl_2})$. In
\cite{Lukpugai} the elliptic algebra is obtained by twisting the
Cartan current. In some sense, we say that the elliptic algebra at
level one governs the structure of the q-Virasoro algebra. It seems
true that it also holds for their higher rank extensions. So
following this approach and the above mentioned expectations, it is
important to obtain the realizations of the elliptic quantum
algebras. In fact, for studying the RSOS model and its higher spin
extension (i.e. the $k$-fusion RSOS model), the representations of
 $U_{q,p}(\widehat{sl_2})$ with any given level $k$ have been
presented in \cite{Konno1} and \cite{CD07}. They are different from
each other. The former can be viewed as the elliptic version of the
parafermionic realization, which is obtained by twisting the
parafermionic realization of the quantum affine algebra
$U_{q}(\widehat{sl_2})$; and the latter is the elliptic deformation
of the Wakimoto realization. The elliptic algebra
$U_{q,p}(\widehat{sl_2})$ is actually the Drinfeld realization of
$B_{q,\lambda}(\widehat{sl_2})$ showed in \cite{JKOS}. Furthermore,
in order to study a higher rank extension of the RSOS model, we
should construct the realizations of $U_{q,p}(\widehat{sl_N})$. It
can be viewed as the Drinfeld realization of the face type elliptic
algebra $B_{q,\lambda}(\widehat{sl_N})$ showed in \cite{JKOS,KK}.
However only in the level-one case, the parafermion realization of
it was given in \cite{KK}. And it can not be extended to the higher
level $k$, although parafermion theory is important in physics
\cite{Nemesch,GepQiu,Gep} and in mathematics \cite{Dong}. The
realizations of \cite{Konno1,KK} are based on the facts that in
$\widehat{su(2)}_k$ case, the parafermions are decoupled from the
Cartan current, while in $\widehat{su(N)}_1$ case, the parafermions
become trivial (i.e.identity operator). In fact, the bosonization of
non-local currents for higher rank and higher-level algebras is a
huge project even in the classical level. So if one wants to deal
with the elliptic quantum algebra of higher rank through
bosonization of the non-local currents, it will not be a practical
way. In this paper we will introduce a new way to construct the free
field representation of the higher rank algebra
$U_{q,p}(\widehat{sl_N})$. It is the higher rank generalization
of the construction in \cite{CD07}. And our construction could be
viewed as a twisted version of the quantum semi-infinite flag
manifolds \cite{FF}.

In free fields approach, there are two necessary ingredients that
one has to discuss: screening currents and vertex operators (VOs).
They all play crucial roles in calculating correlation functions and
investigating the irreducible representations. The screening
currents commute or anti-commute with the currents of
$U_{q,p}(\widehat{sl_N})$ up to a total q-difference of some fields.
And for this algebra, there are two kinds of VOs with distinct
physical applications: the type $I$ VOs and the type $II$ VOs. The
former is a local operator which describes the operation of adding
one lattice site, and the formula of the correlation functions can
be expressed as traces of the product of these operators over
irreducible representation space; while the latter plays the role of
particle creation or annihilation operators. In this paper, we also
construct the free field realization of these two important objects.
In fact, they are all obtained by twisting the corresponding ones of
the quantum affine algebra $U_{q}(\widehat{sl_N})$. In order to do
that, we have to construct the type $II$ VOs of
$U_{q}(\widehat{sl_N})$ which have never been given before. In fact,
even for the classical affine algebra, the type $II$ VOs of it are
unknown.

In this paper, in section $2$ we define the Drinfeld realization of
$U_{q,p}(\widehat{sl_N})$ as a certain tensor product of the quantum
affine algebra $U_{q}(\widehat{sl_N})$ and a Heisenberg algebra,
which is different from the ones given in \cite{JKOS, KK}. With this
definition it is more convenient to construct the free field
representation of $U_{q,p}(\widehat{sl_N})$ with given level $k$.
And in section $3$ we will present the construction in two steps. In
section $4$ a series of screening currents of
$U_{q,p}(\widehat{sl_N})$ are given. In section $5$ the explicit
expressions of the type $II$ VOs of $U_{q}(\widehat{sl_N})$ and
the two types VOs of $U_{q,p}(\widehat{sl_N})$ are presented.

\setcounter{section}{1} \setcounter{equation}{0}

\vskip 1cm
\section{The elliptic quantum algebra $U_{q,p}(\widehat{sl_N})$}
\vskip 0.5cm

There are two types of the elliptic quantum algebras: the face type
and the vertex type. Here we only consider the face type elliptic
algebra $U_{q,p}(\widehat{sl_N})$, which can be viewed as the
Drinfeld realization of the face type elliptic quantum group
$B_{q,\lambda}(\widehat{sl_N})$. Usually, we can also consider it as
the tensor product of the quantum affine algebra
$U_q(\widehat{sl_N})$ and a Heisenberg algebra. In this section, we
will first review the definition of $U_q(\widehat{sl_N})$; then we
will define the elliptic currents of it; lastly, we give the
definition of the elliptic algebra $U_{q,p}(\widehat{sl_N})$.
Throughout this paper, we fix a complex number $q\neq0$, $|q|<1$.

\vskip 0.5cm
\subsection{The quantum affine algebra $U_q(\widehat{sl_N})$}
\vskip 0.5cm

In this subsection, for convenience, we give a review of the
definition of $U_q(\widehat{sl_N})$. We will use the standard symbol
$[n]$:
$$[n]=\frac{q^{n}-q^{-n}}{q-q^{-1}},$$
\no and let $A=(a_{ij})_{1\leq i, j\leq N-1}$ be the Cartan matrix
of $sl_{N}$. The dual Coxeter number of it is denoted by $h^{\vee}$
and $h^{\vee}=N$.

{\defi. $U_q(\widehat{sl_N})$ is the associative algebra over
$\mathbb{C}$ with Drinfeld generators $H_n^i$
$(n\in\mathbb{Z}-\{0\})$, $e_n^{\pm, i}$ $(n\in\mathbb{Z})$,  $h_i$
$(i=1, \cdots, N-1)$ and the central element $c$ satisfying the
following defining relations: \beqa &&[h_i, H_n^j]=0, \ \ \ \ [h_i,
e_n^{\pm, j}]=\pm a_{ij}e_n^{\pm,
j}\\
&&[H_n^i, H_m^j]=\frac{[a_{ij}n][cn]}{n}\delta_{n+m, 0},\\
&&[H_n^i, e_m^{\pm,
j}]=\pm\frac{[a_{ij}n]}{n}q^{\mp\frac{c}{2}n}e_{n+m}^{\pm, j},\\
&&[e_n^{+, i}, e_m^{-,
j}]=\frac{\delta^{ij}}{q-q^{-1}}(q^{\frac{c}{2}(n-m)}\psi_{+,
n+m}^i-q^{-\frac{c}{2}(n-m)}\psi_{-, n+m}^i),\\
&&[e_{n+1}^{\pm, i}, e_m^{\pm, j}]_{q^{\pm a_{ij}}}+[e_{m+1}^{\pm,
j}, e_n^{\pm, i}]_{q^{\pm a_{ij}}}=0,\\
&&[e_n^{\pm, i}, e_m^{\pm, j}]=0 \ \ \  for\ \ \  a_{ij}=0,\\
&&[e_n^{\pm, i}, [e_m^{\pm, i}, e_l^{\pm, j}]_{q^{\mp 1}}]_{q^{\pm
1}}+[e_m^{\pm, i}, [e_n^{\pm, i}, e_l^{\pm, j}]_{q^{\mp 1}}]_{q^{\pm
1}}=0 \ \ \ for\ \ \ a_{ij}=-1, \eeqa

\no where $\psi_{\pm, n}^{i}$ are defined by
$$\sum_{n\in\mathbb{Z}}\psi_{\pm, n}^{i}z^{-n}=q^{\pm
h_i}\exp\Big(\pm(q-q^{-1})\sum_{\pm n>0}H_n^{i}z^{-n}\Big),$$ and
the symbol $[A,B]_{x}$ for $x\in \mathbb{C}$ denotes $AB-xBA.$}

If we introduce the generating functions $\psi_{\pm}^{i}(z)$ and
$e^{\pm, i}(z)$ $(i=1, \cdots, N-1)$ as \beqas
&&\psi_{\pm}^{i}(z)=\sum_{n\in \mathbb{Z}}\psi_{\pm, n}^{i}z^{-n}, \
\ \ \ e^{\pm, i}(z)=\sum_{n\in\mathbb{Z}}e_n^{\pm, i}z^{-n-1},
\eeqas

\no which are called the Drinfeld currents of $U_q(\widehat{sl_N})$.
In terms of them, the above defining relations (2.1)-(2.7) can be
recast as \beqa
&&[\psi_{\pm}^{i}(z), \psi_{\pm}^{j}(w)]=0,\\
&&(z-q^{a_{ij}-c}w)(z-q^{-a_{ij}+c}w)\psi_{+}^{i}(z)\psi_{-}^{j}(w)\nn\\
&&\ \ \ \ \ =(z-q^{a_{ij}+c}w)(z-q^{-a_{ij}-c}w)\psi_{-}^{j}(w)\psi_{+}^{i}(z),\\
&&(z-q^{\pm (a_{ij}-\frac{c}{2})}w)\psi_{+}^{i}(z)e^{\pm, j}(w)=
(q^{\pm a_{ij}}z-q^{\mp \frac{c}{2}}w)e^{\pm, j}(w)\psi_{+}^{i}(z),\\
&&(z-q^{\pm (a_{ij}-\frac{c}{2})}w)e^{\pm, j}(z)\psi_{-}^{i}(w)=
(q^{\pm a_{ij}}z-q^{\mp \frac{c}{2}}w)\psi_{-}^{i}(w)e^{\pm,
j}(z),\\
&&[e^{+, i}(z), e^{-, j}(w)]=\frac{\delta^{ij}}{(q-q^{-1})zw}
\Big(\delta(q^{c}w/z)\psi_{+}^{i}(q^{\frac{c}{2}}w)
-\delta(q^{-c}w/z)\psi_{-}^{i}(q^{-\frac{c}{2}}w)\Big),\\
&&(z-q^{\pm a_{ij}}w)e^{\pm, i}(z)e^{\pm, j}(w)=(q^{\pm
a_{ij}}z-w)e^{\pm, j}(w)e^{\pm, i}(z),\\
&&e^{\pm, i}(z)e^{\pm, j}(w)=e^{\pm, j}(w)e^{\pm, i}(z)\ \ \ for\ \
\ a_{ij}=0,\\
&&e^{\pm, i}(z_1)e^{\pm, i}(z_2)e^{\pm, j}(w)-[2]e^{\pm,
i}(z_1)e^{\pm, j}(w)e^{\pm, i}(z_2)\nn \\
&&\hskip0.5cm +e^{\pm, j}(w)e^{\pm, i}(z_1)e^{\pm,
i}(z_2)+(replacement: z_1\leftrightarrow z_2)=0\ \ \ for\ \ \
a_{ij}=-1, \eeqa

\no where $\dl(x)=\sum _{n\in\mathbb{Z}}x^n$.

\vskip 0.5cm
\subsection{The elliptic algebra $U_{q,p}(\widehat{sl_N})$}
\vskip 0.5cm

The elliptic algebra $U_{q,p}(\widehat{sl_N})$ can be considered as
the tensor product of the elliptic currents of $U_q(\widehat{sl_N})$
and a Heisenberg algebra \cite{JKOS}. We first give the elliptic
currents of $U_q(\widehat{sl_N})$. A pair of parameters $p$ and
$p^*$ will be used:
$$ p=q^{2r}, \ \
p^*=q^{2r^*}=pq^{-2c} \ \ \ \ \ \ \ (r^*=r-c; \ \
r,r^*\in\mathbb{R}_{>0}). $$

Let us define the currents $D_i^{\pm}(z; r, r^*)$ $\in$
$U_q(\widehat{sl_N})$\ \ $(i=1, \cdots, N-1)$ depending on $r$ and
$r^*$ as \beqas &&D_i^{+}(z; r,
r^*)=\exp\Big(\sum_{n>0}\frac{1}{[r^*n]}H_{-n}^{i}q^{(r^{*}+\frac{c}{2})n}z^n\Big),
\\
&&D_i^{-}(z; r,
r^*)=\exp\Big(-\sum_{n>0}\frac{1}{[rn]}H_{n}^{i}q^{(r-\frac{c}{2})n}z^{-n}\Big),
\eeqas

\no which are different from the ones in \cite{JKOS, KK} by a power
of $q$. Using them we can define the `dressed' currents
$\Psi_i^{\pm}(z)$, $e_i(z)$ and $f_i(z)$ $(i=1, \cdots, N-1)$ as:
\beqas &&\Psi_i^{+}(z)=D_i^{+}(q^{\frac{c}{2}}z; r, r^*)\psi
_{+}^i(z)D_i^{-}(q^{-\frac{c}{2}}z; r, r^*), \\
&&\Psi_i^{-}(z)=D_i^{+}(q^{-\frac{c}{2}}z; r, r^*)\psi
_{-}^i(z)D_{i}^{-}(q^{\frac{c}{2}}z; r, r^*), \\
&&e_i(z)=D_i^{+}(z; r, r^*)e^{+, i}(z), \\
&&f_i(z)=e^{-, i}(z)D_i^{-}(z; r, r^*). \eeqas

\no Obviously these currents all depend on the parameter $p$.
Moreover, applying (2.8)-(2.15) we have the following proposition by
direct calculation:

{\prop. The fields $\Psi_i^{\pm}(z)$, $e_i(z)$ and $f_i(z)$ $(i=1,
\cdots, N-1)$ defined above satisfy the following elliptic
commutation relations: \beqa
&&\Psi_i^{\pm}(z)\Psi_j^{\pm}(w)=\frac{\Theta_p(q^{-a_{ij}}\frac{z}{w})
\Theta_{p^*}(q^{a_{ij}}\frac{z}{w})}
{\Theta_p(q^{a_{ij}}\frac{z}{w})\Theta_{p^*}(q^{-a_{ij}}\frac{z}{w})}
\Psi_j^{\pm}(w)\Psi_i^{\pm}(z),\\
&&\Psi_i^{+}(z)\Psi_j^{-}(w)=\frac{\Theta_p(pq^{-a_{ij}-c}\frac{z}{w})
\Theta_{p^*}(p^{*}q^{a_{ij}+c}\frac{z}{w})}
{\Theta_p(pq^{a_{ij}-c}\frac{z}{w})\Theta_{p^*}(p^{*}q^{-a_{ij}+c}\frac{z}{w})}
\Psi_j^{-}(w)\Psi_i^{+}(z),\\
&&\Psi_i^{\pm}(z)e_j(w)=q^{-a_{ij}}\frac{\Theta_{p^*}(q^{\pm
\frac{c}{2}+a_{ij}}\frac{z}{w})}{\Theta_{p^*}(q^{\pm
\frac{c}{2}-a_{ij}}\frac{z}{w})}e_j(w)\Psi_i^{\pm}(z),\\
&&\Psi_i^{\pm}(z)f_i(w)=q^{a_{ij}}\frac{\Theta_p(q^{\mp
\frac{c}{2}-a_{ij}}\frac{z}{w})}{\Theta_p(q^{\mp
\frac{c}{2}+a_{ij}}\frac{z}{w})}f_i(w)\Psi_i^{\pm}(z),\\
&&[e_i(z),
f_j(w)]=\frac{\delta^{ij}}{(q-q^{-1})zw}\Big(\delta(q^{-c}\frac{z}{w})
\Psi_i^{+}(q^{\frac{c}{2}}w)
-\delta(q^{c}\frac{z}{w})\Psi_i^{-}(q^{-\frac{c}{2}}w)\Big),\\
&&e_i(z)e_j(w)=q^{-a_{ij}}\frac{\Theta_{p^*}(q^{a_{ij}}\frac{z}{w})}
{\Theta_{p^*}(q^{-a_{ij}}\frac{z}{w})}
e_j(w)e_i(z),\\
&&f_i(z)f_j(w)=q^{a_{ij}}\frac{\Theta_p(q^{-a_{ij}}\frac{z}{w})}
{\Theta_p(q^{a_{ij}}\frac{z}{w})}f_j(w)f_i(z),\\
&&\frac{(p^{*}q^{2}\frac{z_2}{z_1};
p^*)_\infty}{(p^{*}q^{-2}\frac{z_2}{z_1};
p^*)_\infty}\Bigg\{e_j(w)e_i(z_1)e_i(z_2)-[2]\frac{(p^{*}q\frac{z_1}{w};
p^*)_\infty(p^{*}q^{-1}\frac{w}{z_1};
p^*)_\infty}{(p^{*}q^{-1}\frac{z_1}{w};
p^*)_\infty(p^{*}q\frac{w}{z_1};
p^*)_\infty}e_i(z_1)e_j(w)e_i(z_2)\nn\\
&&\hskip2.2cm +\frac{(p^{*}q\frac{z_1}{w};
p^*)_\infty(p^{*}q^{-1}\frac{w}{z_1};
p^*)_\infty}{(p^{*}q^{-1}\frac{z_1}{w};
p^*)_\infty(p^{*}q\frac{w}{z_1};
p^*)_\infty}\frac{(p^{*}q\frac{z_2}{w};
p^*)_\infty(p^{*}q^{-1}\frac{w}{z_2};
p^*)_\infty}{(p^{*}q^{-1}\frac{z_2}{w};
p^*)_\infty(p^{*}q\frac{w}{z_2};
p^*)_\infty}e_i(z_1)e_i(z_2)e_j(w)\Bigg\}\nn\\
&&\hskip1cm +(replacement: z_1\leftrightarrow z_2)=0\hskip1cm
for\hskip0.3cm
|i-j|\leq 1,\\
&&\frac{(pq^{-2}\frac{z_2}{z_1}; p)_\infty}{(pq^{2}\frac{z_2}{z_1};
p)_\infty}\Bigg\{f_j(w)f_i(z_1)f_i(z_2)-[2]\frac{(pq\frac{w}{z_1};
p)_\infty(pq^{-1}\frac{z_1}{w}; p)_\infty}{(pq^{-1}\frac{w}{z_1};
p)_\infty(pq\frac{z_1}{w};
p)_\infty}f_i(z_1)f_j(w)f_i(z_2)\nn\\
&&\hskip2cm +\frac{(pq\frac{w}{z_1}; p)_\infty(pq^{-1}\frac{z_1}{w};
p)_\infty}{(pq^{-1}\frac{w}{z_1}; p)_\infty(pq\frac{z_1}{w};
p)_\infty}\frac{(pq\frac{w}{z_2}; p)_\infty(pq^{-1}\frac{z_2}{w};
p)_\infty}{(pq^{-1}\frac{w}{z_2}; p)_\infty(pq\frac{z_2}{w};
p)_\infty}f_i(z_1)f_i(z_2)f_j(w)\Bigg\}\nn\\
&&\hskip1cm +(replacement: z_1\leftrightarrow z_2)=0\hskip1cm
for\hskip0.3cm |i-j|\leq 1, \eeqa

\no where we use the elliptic theta function $\Theta_t(z)$ for any
parameter $t=q^{2\nu}\ \ (\nu\in \mathbb{C})$ defined as

$$\Theta_t(z)=(z; t)_{\infty}(tz^{-1}; t)_{\infty}(t; t)_{\infty},$$

\no in which \\

$$(z; t_1,\cdots, t_k)_{\infty}=\prod_{n_1, \cdots,
n_k\geq0}(1-zt_1^{n_1}\cdots t_k^{n_k}).$$}

\no Here these `dressed' currents $\Psi_i^{\pm}(z)$, $e_i(z)$ and
$f_i(z)$ $(i=1, \cdots, N-1)$ are called the elliptic currents of
$U_q(\widehat{sl_N})$ since they obey the above elliptic commutation
relations.

Next, we need a set of Heisenberg algebras generated by {$P_i$,
$Q_i$} $(i=1, \cdots, N-1)$ with \\
$$[P_i, Q_j]=-\frac{a_{ij}}{2},$$ to add nice periodicity
properties to the elliptic exchange relations (2.16)-(2.24). And the
Heisenberg algebras commute with $U_q(\widehat{sl_N})$. For
convenience, the following parametrization will be used in the
following sections: \beqas
&&q=e^{-\pi i/r\tau}, \\
&&p=e^{-2\pi i/\tau}, \ \ \ \ p^*=e^{-2\pi i/\tau^{*}}\\
&&z=q^{2u}=e^{-2\pi iu/r\tau}. \eeqas

\no With them, we can further define the currents $H_{i}^{\pm}(u)$,
$E_{i}(u)$ and $F_{i}(u)$\ \ $(i=1, \cdots, N-1)$ as follows:
 \beqas
&&H_i^{\pm}(u)=\Psi_i^{\pm}(z)e^{2Q_i}q^{\mp h_i}
(q^{\pm(r-\frac{c}{2})}z)^{\frac{(h_{i}+P_{i}-1)}{r}-\frac{(P_{i}-1)}{r^*}},
\\
&&E_i(u)=e_i(z)e^{2Q_i}z^{-\frac{(P_{i}-1)}{r^*}}, \\
&&F_i(u)=f_i(z)z^{\frac{(h_{i}+P_{i}-1)}{r}} \eeqas

\no They are actually the tensor product of elliptic currents
$\Psi_i^{\pm}(z)$, $e_i(z)$ and $f_i(z)$ with the Heisenberg
algebras. And to distinguish them from the elliptic currents, we
call them the total currents. It should be noted that the choice of
the zero-modes in $H_{i}^{\pm}(u)$, $E_{i}(u)$ and $F_{i}(u)$ are
different from the ones given in \cite{JKOS, KK}. Our choice makes
our construction of the free field realization of
$U_{q,p}(\widehat{sl_N})$ more convenience. Now the definition of
$U_{q,p}(\widehat{sl_N})$ can be stated explicitly as:

{\defi. The elliptic algebra $U_{q,p}(\widehat{sl_N})$ is isomorphic
to the associative algebra over $\mathbb{C}$ generated by
$H_{i}^{\pm}(u)$, $E_{i}(u)$ and $F_{i}(u)$\ \ $(i=1, \cdots, N-1)$
with the following defining relations: \beqa
&&H_{i}^{\pm}(u)H_{j}^{\pm}(v)=\frac{\theta_{r}(u-v-\frac{a_{ij}}{2})}
{\theta_{r}(u-v+\frac{a_{ij}}{2})}\frac{\theta_{r^*}(u-v+\frac{a_{ij}}{2})}
{\theta_{r^*}(u-v-\frac{a_{ij}}{2})}H_{j}^{\pm}(v)H_{i}^{\pm}(u),
\\
&&H_{i}^{+}(u)H_{j}^{-}(v)=\frac{\theta_{r}(u-v-\frac{c}{2}-\frac{a_{ij}}{2})}
{\theta_{r}(u-v-\frac{c}{2}+\frac{a_{ij}}{2})}
\frac{\theta_{r^*}(u-v+\frac{c}{2}+\frac{a_{ij}}{2})}
{\theta_{r^*}(u-v+\frac{c}{2}-\frac{a_{ij}}{2})}H_{j}^{-}(v)H_{i}^{+}(u),\\
&&H_{i}^{\pm}(u)E_j(v)=\frac{\theta_{r^*}(u-v\pm\frac{c}{4}+\frac{a_{ij}}{2})}
{\theta_{r^*}(u-v\pm\frac{c}{4}-\frac{a_{ij}}{2})}E_j(v)H_{i}^{\pm}(u),
\\
&&H_{i}^{\pm}(u)F_j(v)=\frac{\theta_{r}(u-v\mp\frac{c}{4}-\frac{a_{ij}}{2})}
{\theta_{r}(u-v\mp\frac{c}{4}+\frac{a_{ij}}{2})}F_j(v)H_{i}^{\pm}(u),\\
&&[E_i(u), F_j(v)]=\frac{\delta_{ij}}{(q-q^{-1})zw}
\Big(\dl(u-v-\frac{c}{2})H^{+}_{i}(u-\frac{c}{4})\nn\\
&&\hskip4.75cm -\dl(u-v+\frac{c}{2})H^{-}_{i}(v-\frac{c}{4})\Big),\\
&&E_i(u)E_j(v)=\frac{\theta_{r^*}(u-v+\frac{a_{ij}}{2})}
{\theta_{r^*}(u-v-\frac{a_{ij}}{2})}E_j(v)E_i(u), \\
&&F_i(u)F_j(v)=\frac{\theta_{r}(u-v-\frac{a_{ij}}{2})}
{\theta_{r}(u-v+\frac{a_{ij}}{2})}F_j(v)F_i(u), \\
&&z_1^{-\frac{2}{r^*}}\frac{(p^{*}q^{2}\frac{z_2}{z_1};
p^*)_\infty}{(p^{*}q^{-2}\frac{z_2}{z_1};
p^*)_\infty}\Bigg\{E_j(u)E_i(u_1)E_i(u_2)\nn\\
&&\hskip2.8cm
-[2](\frac{z}{z_1})^{\frac{a_{ij}}{r^*}}\frac{(p^{*}q^{a_{ij}}\frac{z}{z_1};
p^*)_\infty(p^{*}q^{-a_{ij}}\frac{z_1}{z};
p^*)_\infty}{(p^{*}q^{-a_{ij}}\frac{z}{z_1};
p^*)_\infty(p^{*}q^{a_{ij}}\frac{z_1}{z};
p^*)_\infty}E_i(u_1)E_j(u)E_i(u_2)\nn\\
&&\hskip2.8cm
+(\frac{z}{z_1})^{\frac{a_{ij}}{r^*}}(\frac{z}{z_2})^{\frac{a_{ij}}{r^*}}
\frac{(p^{*}q^{a_{ij}}\frac{z}{z_1};
p^*)_\infty(p^{*}q^{-a_{ij}}\frac{z_1}{z};
p^*)_\infty}{(p^{*}q^{-a_{ij}}\frac{z}{z_1};
p^*)_\infty(p^{*}q^{a_{ij}}\frac{z_1}{z};
p^*)_\infty}\nn\\
&&\hskip3.2cm\times\frac{(p^{*}q^{a_{ij}}\frac{z}{z_2};
p^*)_\infty(p^{*}q^{-a_{ij}}\frac{z_2}{z};
p^*)_\infty}{(p^{*}q^{-a_{ij}}\frac{z}{z_2};
p^*)_\infty(p^{*}q^{a_{ij}}\frac{z_2}{z};
p^*)_\infty}E_i(u_1)E_i(u_2)E_j(u)\Bigg\}\nn\\
&&\hskip1cm +(replacement: z_1\leftrightarrow z_2)=0\hskip1cm
for\hskip0.3cm
|i-j|\leq 1,\\
&&z_1^{\frac{2}{r}}\frac{(pq^{-2}\frac{z_2}{z_1};
p)_\infty}{(pq^{2}\frac{z_2}{z_1};
p)_\infty}\Bigg\{F_j(u)F_i(u_1)F_i(u_2)\nn\\
&&\hskip2.3cm
-[2](\frac{z}{z_1})^{-\frac{a_{ij}}{r}}\frac{(pq^{-a_{ij}}\frac{z}{z_1};
p)_\infty(pq^{a_{ij}}\frac{z_1}{z};
p)_\infty}{(pq^{a_{ij}}\frac{z}{z_1};
p)_\infty(pq^{-a_{ij}}\frac{z_1}{z};
p)_\infty}F_i(u_1)F_j(u)F_i(u_2)\nn\\
&&\hskip2.3cm
+(\frac{z}{z_1})^{-\frac{a_{ij}}{r}}(\frac{z}{z_2})^{-\frac{a_{ij}}{r}}
\frac{(pq^{-a_{ij}}\frac{z}{z_1};
p)_\infty(pq^{a_{ij}}\frac{z_1}{z};
p)_\infty}{(pq^{a_{ij}}\frac{z}{z_1};
p)_\infty(pq^{-a_{ij}}\frac{z_1}{z};
p)_\infty}\nn\\
&&\hskip3cm \times\frac{(pq^{-a_{ij}}\frac{z}{z_2};
p)_\infty(pq^{a_{ij}}\frac{z_2}{z};
p)_\infty}{(pq^{a_{ij}}\frac{z}{z_2};
p)_\infty(pq^{-a_{ij}}\frac{z_2}{z};
p)_\infty}F_i(u_1)F_i(u_2)F_j(u)\Bigg\}\nn\\
&&\hskip1cm +(replacement: z_1\leftrightarrow z_2)=0\hskip1cm
for\hskip0.3cm |i-j|\leq 1, \eeqa

\no where the notation of the Jacobi theta functions
$\theta_{\nu}(u)$ for $\nu\in \mathbb{C}$ are used,
$$\theta_{\nu}(u)=q^{\frac{u^2}{\nu}-u}\frac{\Theta_{q^{2\nu}}(q^{2u})}
{(q^{2\nu}; q^{2\nu})_{\infty}^3}.$$}

\no Note that we have used the parametrization $z=q^{2u}$,
$w=q^{2v}$ and $z_{i}=q^{2u_{i}}$\ \ $(i=1, 2)$ in the above
expressions. In the following, we will use this parametrization
without mentioning them if they are not confused. It is easy to see
that the above relations (2.25)-(2.33) have good periodicity
properties because of the quasi-periodicity property of the Jacobi
theta functions, such as
$$\theta_{r}(u+r)=-\theta_{r}(u), \ \ \ \ \ \ \theta_{r}(u+r\tau)=-e^{-\pi \tau
i-2\pi iu/r}\theta_{r}(u)$$ \no and similar relations hold for
$\theta_{r^*}(u)$ with $r$ replaced by $r^*$.

\setcounter{section}{2} \setcounter{equation}{0} \vskip 1cm
\section{Free field realization of $U_{q,p}(\widehat{sl_N})$}
\vskip 0.5cm

The level $k$ representation of $U_{q,p}(\widehat{sl_N})$ has not
been given before. Although the free field realization of it in
level 1 was given in \cite{KK}, it can not be generalized to the
higher level case. In this section, by using a new method, we will
construct a free boson realization of $U_{q,p}(\widehat{sl_N})$ with
given level $k$. This method has been used to construct a free field
realization of $U_{q,p}(\widehat{sl_2})_k$ in \cite{CD07}. Here we
will show that it can be generalized to the higher rank case. The
method is to twist the level $k$ Wakimoto realization of
$U_q(\widehat{sl_N})$ by constructing some `twising' currents. We
will first fix some conventions and review the Wakimoto realization
of $U_q(\widehat{sl_N})$ in \cite{AOS}; then we will give our
construction in two steps: the first one is the bosonization of the
elliptic currents of $U_q(\widehat{sl_N})$; and the second one is
the free boson realization of the total currents.

\vskip 0.5cm
\subsection{Notations}
\vskip 0.5cm

We introduce a quantum Heisenberg algebra $\mathscr{H}_{q, k}$ with
the generators: $a_n^i, p_{a}^{i}, q_{a}^{i}$ for $1\leq i\leq N-1$;
$b_{n}^{ij}, p_{b}^{ij}, q_{b}^{ij}$ and $c_{n}^{ij}, p_{c}^{ij},
q_{c}^{ij}$ for $1\leq i<j\leq N$, where $n\in \mathbb{Z}_{\neq 0}$,
and the defining relations are as follows: \beqas
&&[a_n^i,a_m^j]=\frac{[(k+h^{\vee})n][a_{ij}n]}{n}\dl_{n+m, 0},\ \ \
[p_{a}^{i},
q_{a}^j]=a_{ij}(k+h^{\vee}),\\
&&[b_{n}^{ij},
b_{m}^{i'j'}]=-\frac{[n]^2}{n}\delta^{ii'}\delta^{jj'}
\dl_{n+m, 0},\ \ \ [p_{b}^{ij}, q_{b}^{i'j'}]=-\delta^{ii'}\delta^{jj'},\\
&&[c_{n}^{ij},
c_{m}^{i'j'}]=\frac{[n]^2}{n}\delta^{ii'}\delta^{jj'}\dl_{n+m, 0},\
\ \ [p_{c}^{ij}, q_{c}^{i'j'}]=\delta^{ii'}\delta^{jj'} \eeqas and
the others vanish. Using them, we set the generating functions
$a^{i}(z; \alpha)$ for $\alpha \in \mathbb{C}$ and $a_{\pm}^{i}(z)$
$(1\leq i\leq N-1)$ by: \beqas
&&a^{i}(z; \alpha)=-\sum_{n\neq0}\frac{a_n^i}{[n]}q^{-\alpha |n|}z^{-n}+q_{a}^{i}+p_{a}^{i}\ln z,\\
&&a_{\pm}^{i}(z)=\pm\Big((q-q^{-1})\sum_{n>0}a_{\pm n}^{i}z^{\mp
n}+p_{a}^{i}\ln q\Big) \eeqas \no and $a^{i}(z; 0) \equiv a^{i}(z)$
for simplicity. Similarly, the generating functions $b^{i, j}(z;
\alpha)$, $b_{\pm}^{i, j}(z)$ and $c^{i, j}(z; \alpha)$,
$c_{\pm}^{i, j}(z)$ for $1\leq i<j\leq N$ can also be given. These
generating functions can be viewed as some free bosonic fields, if
we consider the generators of $\mathscr{H}_{q, k}$ as the modes of
$N^{2}-1$ free bosons: $a^i$ $(1\leq i\leq N-1)$, $b^{ij}$ and
$c^{ij}$ $(1\leq i<j\leq N)$. We also define the completion
$\widehat{\mathscr{H}}_{q, k}$ of $\mathscr{H}_{q, k}$ as:

$$\widehat{\mathscr{H}}_{q, k}=\lim_{\leftarrow}\mathscr{H}_{q, k}/I_{n},\ \  n> 0,$$

\no where $I_{n}$ is the left ideal of $\mathscr{H}_{q, k}$
generated by all the polynomials in \{$a_m^i (1\leq i\leq N-1)$,
$b_{m}^{ij}$ and $c_{m}^{ij} (1\leq i<j\leq N)$: $m> 0$\} of degree
greater than or equal to $n$ \   (here we set
deg$(a_m^i)$=deg$(b_{m}^{ij})$=\\deg$(c_{m}^{ij})$=m). Normal order
prescription $: \ :$ is set by moving $a_n^i (n>0)$ and $p_{a}^{i}$
to the right, while moving $a_n^i (n<0)$ and $q_{a}^{i}$ to the
left. For example,
$$: \exp(a^{i}(z)) :=\exp\Big(-\sum_{n<0}\frac{a_n^i}{[n]}z^{-n}\Big)
e^{q_{a}^{i}}z^{p_{a}^{i}}
\exp\Big(-\sum_{n>0}\frac{a_n^i}{[n]}z^{-n}\Big).$$

In terms of the above free bosonic fields, we can define a
homomorphism $h_{q, k}$ from the algebra $U_q(\widehat{sl_N})$
to $\widehat{\mathscr{H}}_{q, k}$. It is defined on the generators
by: \beqa &&h_{q,
k}\Big(\psi_{\pm}^{i}(z)\Big)=:\exp\Big(\sum_{j=1}^{i}(b_{\pm}^{j,
i+1}(q^{\pm(\frac{k}{2}+j-1)}z)-b_{\pm}^{j,
i}(q^{\pm(\frac{k}{2}+j)}z))
+a_{\pm}^{i}(q^{\pm{\frac{h^{\vee}}{2}}}z)\nn\\
&&\hskip2.5cm +\sum_{j=i+1}^{N}(b_{\pm}^{i,
j}(q^{\pm(\frac{k}{2}+j)}z)-b_{\pm}^{i+1,
j}(q^{\pm(\frac{k}{2}+j-1)}z))\Big):,\\
&&h_{q, k}\Big(e^{+,
i}(z)\Big)=\frac{-1}{(q-q^{-1})z}\sum_{j=1}^{i}:\exp\Big((b+c)^{j,
i}(q^{j-1}z)\Big)\nn \\
&&\hskip3.4cm \times\Big(\exp(b_{+}^{j, i+1}(q^{j-1}z)-(b+c)^{j,
i+1}(q^{j}z))\nn
\\
&&\hskip3.6cm -\exp(b_{-}^{j, i+1}(q^{j-1}z)-(b+c)^{j, i+1}(q^{j-2}z))\Big)\nn \\
&&\hskip3.4cm \times\exp\Big(\sum_{l=1}^{j-1}(b_{+}^{l,
i+1}(q^{l-1}z)-b_{+}^{l,
i}(q^{l}z))\Big):,\\
&&h_{q, k}\Big(e^{-,
i}(z)\Big)=\frac{-1}{(q-q^{-1})z}\Bigg(\sum_{j=1}^{i-1}:\exp\Big((b+c)^{j,
i+1}(q^{-(k+j)}z)\Big)\nn \\
&&\hskip3.8cm \times\Big(\exp(-b_{-}^{j, i}(q^{-(k+j)}z)-(b+c)^{j,
i}(q^{-(k+j-1)}z))\nn \\
&&\hskip4.0cm -\exp(-b_{+}^{j, i}(q^{-(k+j)}z)-(b+c)^{j,
i}(q^{-(k+j+1)}z))\Big)\nn
\\
&&\hskip3.8cm\ \times\exp\Big(\sum_{l=j+1}^{i}(b_{-}^{l,
i+1}(q^{-(k+l-1)}z)-b_{-}^{l, i}(q^{-(k+l)}z))\nn \\
&&\hskip4.0cm
+a_{-}^{i}(q^{-\frac{k+h^{\vee}}{2}}z)+\sum_{l=i+1}^{N}(b_{-}^{i,
l}(q^{-(k+l)}z)-b_{-}^{i+1, l}(q^{-(k+l-1)}z))\Big):\nn
\\
&&\hskip3.2cm +:\exp\Big((b+c)^{i, i+1}(q^{-(k+i)}z)\Big)\nn \\
&&\hskip3.6cm
\times\exp\Big(a_{-}^{i}(q^{-\frac{k+h^{\vee}}{2}}z)+\sum_{l=i+1}^{N}(b_{-}^{i,
l}(q^{-(k+l)}z)-b_{-}^{i+1, l}(q^{-(k+l-1)}z))\Big):\nn \\
&&\hskip3.2cm -:\exp\Big((b+c)^{i, i+1}(q^{k+i}z)\Big)\nn \\
&&\hskip3.6cm
\times\exp\Big(a_{+}^{i}(q^{\frac{k+h^{\vee}}{2}}z)+\sum_{l=i+1}^{N}(b_{+}^{i,
l}(q^{k+l}z)-b_{+}^{i+1, l}(q^{k+l-1}z))\Big):\nn \\
&&\hskip3.2cm -\sum_{j=i+2}^{N}:\exp\Big((b+c)^{i, j}(q^{k+j-1}z)\Big)\nn \\
&&\hskip3.8cm \times\Big(\exp(b_{+}^{i+1, j}(q^{k+j-1}z)-(b+c)^{i+1,
j}(q^{k+j}z))\nn \\
&&\hskip4.0cm -\exp(b_{-}^{i+1, j}(q^{k+j-1}z)-(b+c)^{i+1,
j}(q^{k+j-2}z))\Big)\nn
\\
&&\hskip3.8cm
\times\exp\Big(a_{+}^{i}(q^{\frac{k+h^{\vee}}{2}}z)+\sum_{l=j}^{N}(b_{+}^{i,
l}(q^{k+l}z)-b_{+}^{i+1, l}(q^{k+l-1}z))\Big):\Bigg). \eeqa

\no Then we have the next proposition followed from \cite{AOS}:

{\prop. $h_{q, k}\Big(\psi_{\pm}^{i}(z)\Big)$ and $h_{q,
k}\Big(e^{\pm, i}(z)\Big)$ $(i=1, \cdots N-1)$ with $k=c$ satisfy
the commutation relations (2.8)-(2.15).}

As a result, when $k\neq -h^{\vee}$, this homomorphism $h_{q, k}$
gives the Wakimoto realization of the quantum affine algebra
$U_q(\widehat{sl_N})$ with $k=c$. In the following subsections, we
will construct the free field realization of the elliptic algebra
$U_{q,p}(\widehat{sl_N})$ by twisting this realization.

\vskip 0.5cm
\subsection{Bosonization of elliptic currents}
\vskip 0.5cm

In this subsection, we show the first step of our construction:
giving the bosonization of the elliptic currents $\Psi_i^{\pm}(z)$,
$e_i(z)$ and $f_i(z)$\ \ $(i=1, \cdots N-1)$ of
$U_q(\widehat{sl_N})$. For brevity, in what follows we will use the
same notations for the elements of $U_q(\widehat{sl_N})$ and their
images in the completion of $\mathscr{H}_{q, k}$. Here we need to
introduce some new currents $D_i^{\pm}(z; r, r^*)$\ \ $(i=1, \cdots
N-1)$ depending on parameters $r$ and $r^*$ as: \beqas &&D_i^{+}(z;
r,
r^*)=\exp\Bigg\{\sum_{n>0}\frac{1}{[r^{*}n]}\Big(\sum_{j=1}^{i}(b_{-n}^{j,
i+1}q^{-(k+j-1)n}-b_{-n}^{j,
i}q^{-(k+j)n})+a_{-n}^{i}q^{-\frac{k+h^{\vee}}{2}n}\\
&&\hskip3.2cm +\sum_{j=i+1}^{N}(b_{-n}^{i,
j}q^{-(k+j)n}-b_{-n}^{i+1,
j}q^{-(k+j-1)n})\Big)q^{rn}z^n\Bigg\},\\
&&D_i^{-}(z; r,
r^*)=\exp\Bigg\{-\sum_{n>0}\frac{1}{[rn]}\Big(\sum_{j=1}^{i}(b_{n}^{j,
i+1}q^{-(k+j-1)n}-b_{n}^{j,
i}q^{-(k+j)n})+a_{n}^{i}q^{-\frac{k+h^{\vee}}{2}n}\\
&&\hskip3.2cm +\sum_{j=i+1}^{N}(b_{n}^{i, j}q^{-(k+j)n}-b_{n}^{i+1,
j}q^{-(k+j-1)n})\Big)q^{rn}z^{-n}\Bigg\}, \eeqas

\no which are nominated as twisting currents; then we have the
following lemma:

{\lem. The currents $D_i^{\pm}(z; r, r^*)$\ \ $(i=1, \cdots N-1)$
and the fields in Eqns.(3.1)-(3.3) satisfy the following commutation
relations: \beqa &&D_i^{+}(z; r, r^*)D_j^{-}(w; r,
r^*)=\frac{(pq^{-a_{ij}-k}\frac{z}{w};
p)_\infty}{(pq^{a_{ij}-k}\frac{z}{w};
p)_\infty}\frac{(p^*q^{a_{ij}+k}\frac{z}{w};
p^*)_\infty}{(p^*q^{-a_{ij}+k}\frac{z}{w}; p^*)_\infty}D_j^{-}(w; r,
r^*)D_i^{+}(z; r, r^*), \\
&&D_i^{\pm}(z; r, r^*)D_j^{\pm}(w; r,
r^*)=D_j^{\pm}(w; r, r^*)D_i^{\pm}(z; r, r^*),\\
&&D_i^{+}(z; r,
r^*)\psi_{+}^{j}(w)=\frac{(pq^{-a_{ij}-\frac{k}{2}}\frac{z}{w};
p^*)_\infty}{(pq^{a_{ij}-\frac{k}{2}}\frac{z}{w};
p^*)_\infty}\frac{(p^{*}q^{a_{ij}-\frac{k}{2}}\frac{z}{w};
p^*)_\infty}{(p^{*}q^{-a_{ij}-\frac{k}{2}}\frac{z}{w};
p^*)_\infty}\psi_{+}^{j}(w)D_i^{+}(z; r, r^*),\\
&&D_i^{+}(z; r, r^*)\psi_{-}^{j}(w)=\psi_{-}^{j}(w)D_i^{+}(z; r,
r^*), \\
&&D_i^{+}(z; r, r^*)e^{+, j}(w)=\frac{(p^*q^{a_{ij}}\frac{z}{w};
p^*)_\infty}{(p^*q^{-a_{ij}}\frac{z}{w}; p^*)_\infty}e^{+,
j}(w)D_i^{+}(z; r, r^*),\\
&&D_i^{+}(z; r, r^*)e^{-, j}(w)=\frac{(p^*q^{-a_{ij}+k}\frac{z}{w};
p^*)_\infty}{(p^*q^{a_{ij}+k}\frac{z}{w}; p^*)_\infty}e^{-,
j}(w)D_i^{+}(z; r, r^*),\\
&&D_i^{-}(z; r, r^*)\psi_{+}^{j}(w)=\psi_{+}^{j}(w)D_i^{-}(z; r,
r^*), \\
&&D_i^{-}(z; r,
r^*)\psi_{-}^{j}(w)=\frac{(pq^{a_{ij}+\frac{k}{2}}\frac{w}{z};
p)_\infty}{(pq^{-a_{ij}+\frac{k}{2}}\frac{w}{z};
p)_\infty}\frac{(p^{*}q^{-a_{ij}+\frac{k}{2}}\frac{w}{z};
p)_\infty}{(p^{*}q^{a_{ij}+\frac{k}{2}}\frac{w}{z};
p)_\infty}\psi_{-}^{j}(w)D_i^{-}(z; r, r^*),\\
&&D_i^{-}(z; r, r^*)e^{+, j}(w)=\frac{(pq^{-a_{ij}-k}\frac{w}{z};
p)_\infty}{(pq^{a_{ij}-k}\frac{w}{z}; p)_\infty}e^{+,
j}(w)D_i^{-}(z; r, r^*),\\
&&D_i^{-}(z; r, r^*)e^{-, j}(w)=\frac{(pq^{a_{ij}}\frac{w}{z};
p)_\infty}{(pq^{-a_{ij}}\frac{w}{z}; p)_\infty}e^{-, j}(w)D_i^{-}(z;
r, r^*). \eeqa

\no Proof:} A straightforward but lengthy operator product
expansions (OPE) calculation verifies this lemma. Here, we only take
the first one as an example. It is obvious to see that
$$D_i^{+}(z; r, r^*)D_j^{-}(w; r,
r^*)=:D_i^{+}(z; r, r^*)D_j^{-}(w; r, r^*):;$$ \no and using the
following formulas: \beqas &&e^{A}e^{B}=e^{[A, B]}e^{B}e^{A}, \ \ \
if\ [A, B]\ \
commute\ with\ A \ and\ B;\\
&&\exp\Big(-\sum_{n>0}\frac{x^{n}}{n}\Big)=1-x;\\
&&(1-x)^{-1}=\sum_{n\geq 0}x^{n},\eeqas \no we can prove the
following relations for three cases: $j=i$, $| j-i |=1$ and $| j-i
|\geq 2$:
$$D_j^{-}(w; r, r^*)D_i^{+}(z; r, r^*)=\frac{(pq^{a_{ij}-k}\frac{z}{w};
p)_\infty}{(pq^{-a_{ij}-k}\frac{z}{w};
p)_\infty}\frac{(p^*q^{-a_{ij}+k}\frac{z}{w};
p^*)_\infty}{(p^*q^{a_{ij}+k}\frac{z}{w}; p^*)_\infty}:D_j^{-}(w; r,
r^*)D_i^{+}(z; r, r^*):,$$ \no then we obtain (3.4) since
$$:D_i^{+}(z; r, r^*)D_j^{-}(w; r, r^*):=:D_j^{-}(w; r,
r^*)D_i^{+}(z; r, r^*):.$$

\no The others can be proved similarly.\ \ \ \ \ \ \ \ \ \ \ \ \ \ \
\ $\Box$

Now twisting the free boson realization (3.1)-(3.3) of
$U_q(\widehat{sl_N})$ with $D_i^{\pm}(z; r, r^*)$, we have free
bosonic fields $\Psi_i^{\pm}(z)$, $e_i(z)$ and $f_i(z)$\ \  $(i=1,
\cdots N-1)$ given by \beqa
&&\Psi_i^{+}(z)=D_i^{+}(q^{\frac{k}{2}}z; r, r^*)
\psi_{+}^i(z)D_i^{-}(q^{-\frac{k}{2}}z; r, r^*),\\
&&\Psi_i^{-}(z)=D_i^{+}(q^{-\frac{k}{2}}z; r,
r^*)\psi_{-}^i(z)D_i^{-}(q^{\frac{k}{2}}z; r, r^*),\\
&&e_i(z)=D_i^{+}(z; r, r^*)e^{+, i}(z),\\
&&f_i(z)=e^{-, i}(z)D_i^{-}(z; r, r^*); \eeqa

\no then applying the above Lemma $1$ and the Proposition $2$, we
can obtain the following theorem:

{\thm. The fields (3.14)-(3.17) with $k=c$ satisfy the elliptic
commutation relations (2.16)-(2.24) in Proposition 1.

\no Proof:} For example, we just prove (2.16). By (3.14),
$$\Psi_i^{+}(z)\Psi_j^{+}(w)=D_i^{+}(q^{\frac{k}{2}}z; r,
r^*)\psi_{+}^i(z)D_i^{-}(q^{-\frac{k}{2}}z; r, r^*) \times
D_j^{+}(q^{\frac{k}{2}}w; r,
r^*)\psi_{+}^j(w)D_j^{-}(q^{-\frac{k}{2}}w; r, r^*),$$

\no and the Proposition $2$ and (3.4)-(3.6) in Lemma $1$, \beqas
&&\psi_{+}^i(z)\psi_{+}^j(w)=\psi_{+}^j(w)\psi_{+}^i(z),\\
&&D_i^{+}(q^{\frac{k}{2}}z; r, r^*)D_j^{-}(q^{-\frac{k}{2}}w; r,
r^*)=\frac{(pq^{-a_{ij}}z/w; p)_\infty(p^*q^{a_{ij}+2k}z/w;
p^*)_\infty}{(pq^{a_{ij}}z/w; p)_\infty(p^*q^{-a_{ij}+2k}z/w;
p^*)_\infty}D_j^{-}(q^{-\frac{k}{2}}w; r,
r^*)D_i^{+}(q^{\frac{k}{2}}z; r, r^*),\\
&&D_i^{+}(q^{\frac{k}{2}}z; r,
r^*)\psi_{+}^{j}(w)=\frac{(pq^{-a_{ij}}z/w;
p^*)_\infty(p^{*}q^{a_{ij}}z/w; p^*)_\infty}{(pq^{a_{ij}}z/w;
p^*)_\infty(p^{*}q^{-a_{ij}}z/w;
p^*)_\infty}\psi_{+}^{j}(w)D_i^{+}(q^{\frac{k}{2}}z; r, r^*), \eeqas

\no we get
$$\Psi_i^{+}(z)\Psi_j^{+}(w)=\frac{(pq^{-a_{ij}}z/w;
p)_\infty(pq^{a_{ij}}w/z; p)_\infty(p^*q^{a_{ij}}z/w;
p^*)_\infty(p^*q^{-a_{ij}}w/z; p^*)_\infty}{(pq^{a_{ij}}z/w;
p)_\infty(pq^{-a_{ij}}w/z; p)_\infty(p^*q^{-a_{ij}}z/w;
p^*)_\infty(p^*q^{a_{ij}}w/z;
p^*)_\infty}\Psi_j^{+}(w)\Psi_i^{+}(z);$$ \no moreover, since the
following identity holds:
$$\frac{(pq^{-a_{ij}}z/w;
p)_\infty(p^*q^{a_{ij}}z/w; p^*)_\infty}{(pq^{a_{ij}}z/w;
p)_\infty(p^*q^{-a_{ij}}z/w; p^*)_\infty}=\frac{(q^{-a_{ij}}z/w;
p)_\infty(q^{a_{ij}}z/w; p^*)_\infty}{(q^{a_{ij}}z/w;
p)_\infty(q^{-a_{ij}}z/w; p^*)_\infty}$$

\no the commutation relation (2.16) is obtained: \beqas
&&\Psi_i^{+}(z)\Psi_j^{+}(w)=\frac{(q^{-a_{ij}}z/w;
p)_\infty(pq^{a_{ij}}w/z; p)_\infty(q^{a_{ij}}z/w;
p^*)_\infty(p^*q^{-a_{ij}}w/z; p^*)_\infty}{(q^{a_{ij}}z/w;
p)_\infty(pq^{-a_{ij}}w/z; p)_\infty(q^{-a_{ij}}z/w;
p^*)_\infty(p^*q^{a_{ij}}w/z;
p^*)_\infty}\Psi_j^{+}(w)\Psi_i^{+}(z)\\
&&\hskip2.15cm =\frac{\Theta_p(q^{-a_{ij}}z/w)
\Theta_{p^*}(q^{a_{ij}}z/w)}
{\Theta_p(q^{a_{ij}}z/w)\Theta_{p^*}(q^{-a_{ij}}z/w)}\Psi_j^{+}(w)\Psi_i^{+}(z).
\eeqas \no The commutation relations (2.17)-(2.24) can be verified
in the same way.\ \ \ \ \ \ $\Box$

{\cor. $\Psi_i^{\pm}(z)$, $e_i(z)$ and $f_i(z)$\ \  $(i=1, \cdots
N-1)$ defined above realize the elliptic currents of
$U_q(\widehat{sl_N})$ with level $k=c$.}

Actually, in the $p\rightarrow 0$ limit, $\Psi_i^{\pm}(z)$, $e_i(z)$
and $f_i(z)$\ \  $(i=1, \cdots N-1)$ give a new free field
representation of $U_q(\widehat{sl_N})$, which is different from
the one in subsection 3.1. More precisely, as $p\rightarrow 0$\ (or
$r\rightarrow \infty$): \beqas &&\Psi_i^{+}(z)\rightarrow
\Big(\psi_{-}^{i}(q^{k}z)\Big)^{-1}, \ \ \ \ \
\Psi_i^{-}(z)\rightarrow \Big(\psi_{+}^{i}(q^{k}z)\Big)^{-1},\\
&&e_i(z)\rightarrow q^{-h_i}\Big(\psi_{-}^{i}(q^{k/2}z)\Big)^{-1}e^{+, i}(z),\\
&&f_i(z)\rightarrow e^{-,
i}(z)q^{h_i}\Big(\psi_{+}^{i}(q^{k/2}z)\Big)^{-1}, \eeqas

\no here \beqa &&h_i=\sum_{l=1}^{i}(p_{b}^{l, i+1}-p_{b}^{l,
i})+p_{a}^{i}+\sum_{l=i+1}^{N}(p_{b}^{i, l}-p_{b}^{i+1, l}),\eeqa
\no which has a lot of useful properties. And we will discuss them
and apply them in the following sections.

\vskip 0.5cm
\subsection{Free field realization of $U_{q,p}(\widehat{sl_N})$}
\vskip 0.5cm

The second step of the construction is presented in this subsection.
We will construct the free boson realization of the total currents.
In order to do that, we need to introduce a Heisenberg algebra
$\mathscr{H}$ generated by $\hat{p}_i$ and $\hat{q}_i$ ($1\leq i\leq
N-1$) such that
$$[\hat{q}_i, \hat{p}_j]=\frac{a_{ij}}{2}$$
and they commute with $a^i$ $(1\leq i\leq N-1)$, $b^{ij}$ and
$c^{ij}$ $(1\leq i<j\leq N)$.

With them we define the fields $H_i^{\pm}(u)$, $E_i(u)$ and
$F_i(u)$\ \ $(i=1, \cdots N-1)$ by \beqa
&&H_i^{\pm}(u)=\Psi_i^{\pm}(z)e^{2\hat{q}_i}q^{\mp h_i}(q^{\pm
(r-\frac{k}{2})}z)^{\frac{(\hat{p}_i+h_i-1)}{r}-\frac{(\hat{p}_i-1)}{r^*}},
\\
&&E_i(u)=e_i(z)e^{2\hat{q}_i}z^{-\frac{(\hat{p}_i-1)}{r^*}},\\
&&F_i(u)=f_i(z)z^{\frac{(h_i+\hat{p}_i-1)}{r}},\eeqa

\no where $h_i$ is given by (3.18). Then (3.19)-(3.21) define a
homomorphism from $U_{q,p}(\widehat{sl_N})$ to
$\widehat{\mathscr{H}}_{q, k}\otimes \mathscr{H}$. Here we have the
following lemma about the $h_i's$: {\lem. The following commutation
relations between $h_i$ and $\Psi_i^{\pm}(z)$, $e_i(z)$ and $f_i(z)$
in Eqns. (3.14)-(3.17)\ \ for $i=1, \cdots N-1$ hold:
\beqas &&[h_i, \Psi_j^{\pm}(z)]=0,\\
&&[h_i, e_j(z)]=a_{ij}e_j(z),\\
&&[h_i, f_j(z)]=-a_{ij}f_j(z). \eeqas}

\no This lemma can be easily verified by using the Hausdorff
formula; and they are the useful properties that the $h_i's$
possess, which we mentioned at the end of the above subsection. Then
we can apply this lemma and the Theorem $1$ to obtain the results
below:

{\thm. The fields in Eqns.(3.19)-(3.21) with $k=c$ obey the
commutation relations given by (2.25)-(2.33).}

{\cor. $H_i^{\pm}(u)$, $E_i(u)$ and $F_i(u)$\ \ $(i=1, \cdots N-1)$
defined above give the free boson realization of
$U_{q,p}(\widehat{sl_N})$ with given level $k=c$.}

\vskip 1cm
\section{Screening currents}
\vskip 0.5cm

In free fields approach, one has to discuss two necessary
ingredients: screening currents and vertex operators. We will only
consider the screening currents of the elliptic quantum algebra
$U_{q,p}(\widehat{sl_N})$ in this section. In 2D CFT, screening
current is a primary field of the energy-momentum tensor with
conformal weight $1$, and its integration gives the screening
charge. It has the property that it commutes with the currents
modulo a total differential of certain field. This property ensures
that the screening charge may be inserted in the correlators by
changing their conformal charges without affecting their conformal
properties. In this section, using the bosons $a^i$ $(1\leq i\leq
N-1)$, $b^{ij}$ and $c^{ij}$ $(1\leq i<j\leq N)$, we will construct
a series of screening currents $S^i(z)$ $(1\leq i\leq N-1)$ of
$U_{q,p}(\widehat{sl_N})$. These currents commute with the currents
modulo a total q-difference of some fields, so they could be
regarded as a quantum deformation of the screening currents in 2D
CFT.

We denote a sort of q-difference operator with a parameter $n\in
\mathbb{Z}_{>0}$ by
$${}_n\partial_{z}X(z)\equiv\frac{X(q^{n}z)-X(q^{-n}z)}{(q-q^{-1})z},$$

\no which is called a total q-difference of a function $X(z)$.
Moreover, to eliminate the total q-difference, one can define the
Jackson integral as
$$\int_{0}^{s\infty}X(z)d_{p}z\equiv s(1-p)\sum_{n\in
\mathbb{Z}}X(sp^{n})p^{n}$$

\no for a scalar $s\in \mathbb{C} \backslash \{0\}$ and a complex
number $p$ such that $\mid p \mid<1$. So that,
$$\int_{0}^{s\infty}({}_n\partial_{z}X(z))d_{p}z=0,$$

\no if it is convergent and we take $p=q^{2n}$.

For simplicity, we set boson fields $A_{\pm}^{i}(L_1, \cdots, L_s;
M_1, \cdots, M_{s+1}|z; \alpha)$\ \ $(i=1, \cdots, N-1)$ for $\alpha
\in \mathbb{C}$ with parameters $L_i$ and $M_j$ $(i, j \in
\mathbb{N})$ as follows: \beqas &&A_{+}^{i}(L_1, \cdots, L_s; M_1,
\cdots, M_{s+1}|z;
\alpha)=\sum_{n>0}\frac{[L_{1}n]\cdots[L_{s}n]}{[M_{1}n]
\cdots[M_{s+1}n]}a_{n}^{i}(q^{\alpha}z)^{-n},\\
&&A_{-}^{i}(L_1, \cdots, L_s; M_1, \cdots, M_{s+1}|z;
\alpha)=\sum_{n>0}\frac{[L_{1}n]\cdots[L_{s}n]}{[M_{1}n]
\cdots[M_{s+1}n]}a_{-n}^{i}(q^{\alpha}z)^{n}, \eeqas

\no then in terms of these boson fields and the ones introduced
before, we express the screening currents $S^{i}(z)$\ \ $(i=1,
\cdots, N-1)$ as: \beqas
&&S^{i}(z)=\frac{-1}{(q-q^{-1})z}:\exp\Big\{A_{-}^{i}(-(k+h^{\vee})|z;
-\frac{k+h^{\vee}}{2})+A_{+}^{i}(k+h^{\vee}|z; \frac{k+h^{\vee}}{2}) \\
&&\hskip4.25cm  -\frac{1}{k+h^{\vee}}(q_{a}^{i}+p_{a}^{i}\ln z)\Big\}:\\
&&\hskip3cm \times\Bigg\{\sum_{j=i+1}^{N}:\exp\Big((b+c)^{i+1,
j}(q^{N-j}z)\Big)\\
&&\hskip4cm \times\Big(\exp(-b_{-}^{i, j}(q^{N-j}z)-(b+c)^{i,
j}(q^{N-j+1}z))\\
&&\hskip4.25cm -\exp(-b_{+}^{i, j}(q^{N-j}z)-(b+c)^{i,
j}(q^{N-j-1}z))\Big)\\
&&\hskip4.25cm \times\exp \Big(\sum_{l=j+1}^{N}(b_{-}^{i+1,
l}(q^{N-l+1}z)-b_{-}^{i, l}(q^{N-l}z))\Big):\Bigg\}, \eeqas

\no and they possess the following properties:

{\thm: $S^{i}(z)$ and the currents $H_i^{\pm}(z)$, $E_i(z)$ and
$F_i(z)$\ \ $(i=1, \cdots N-1)$ in Eqns.(3.19)-(3.21) satisfy the
relations:
\beqas &&H_i^{\pm}(z)S^{j}(w)=S^{j}(w)H_i^{\pm}(z)=O(1),\\
&&E_i(z)S^{j}(w)=S^{j}(w)E_i(z)=O(1),\\
&&F_i(z)S^{j}(w)=S^{j}(w)F_i(z)=\delta^{ij}{}_{(k+h^{\vee})}
\partial_{w}\Big[\frac{1}{z-w}\tilde{S}^{i}(z)\Big]+O(1),\\
&&S^{i}(z)S^{j}(w)=\frac{\theta_{k+h^{\vee}}(u-v+\frac{a_{ij}}{2})}
{\theta_{k+h^{\vee}}(u-v-\frac{a_{ij}}{2})}S^{j}(w)S^{i}(z), \eeqas

\no where the symbol $O(1)$ means regularity and $\tilde{S}^{i}(z)\
\ (i=1, \cdots, N-1)$ are given by: \beqas
&&\tilde{S}^{i}(z)=:\exp\Big\{A_{-}^{i}\Big(-(k+h^{\vee})|z;
\frac{k+h^{\vee}}{2}\Big)+A_{+}^{i}\Big(k+h^{\vee}|z; -\frac{k+h^{\vee}}{2}\Big)\\
&&\hskip 2cm -\frac{1}{k+h^{\vee}}(q_{a}^{i}+p_{a}^{i}\ln
z)\Big\}D_{i}^{-}(z; r, r^*)z^{\frac{(h_{i}+\hat{p}_{i}-1)}{r}}:.
\eeqas}

It is obvious to note that since these screening currents do not
contain the parameter $p$, they are also the screening currents of
the quantum affine algebra $U_q(\widehat{sl_N})$. As a result, the
above theorem can be easily proved by applying the results in
\cite{AOS}. Once we give the explicit expressions of the screening
currents, we can calculate the cohomology and study the
irreducibility of modules of the algebra, which we will discuss
separately in the future.

\vskip 1cm
\section{Vertex operators}
\vskip 0.5cm

In this section, except for the screening currents, we will study
the other important object that one has to discuss in free fields
approach: vertex operators (VOs) of $U_{q,p}(\widehat{sl_N})$. In
WZW model, the primary fields could be realized as the highest
weight representation of Kac-Moody algebra, which are commonly
called as vertex operators (VOs)or intertwining operators. For
quantum affine algebra, in \cite{FR}, the authors defined q-deformed
VOs as certain intertwining operators, which could be regarded as
the quantum counterpart of the primary field in 2D CFT. They play
crucial roles in calculating correlation functions. Following this
approach, in this section we will construct the free field
realization of the VOs of $U_{q,p}(\widehat{sl_N})$. There are two
types of them: the type $I$ VOs and the type $II$ VOs. They can all
be viewed as the elliptic analogs of the primary fields.  The
explicit expressions of them are obtained by twisting the
corresponding ones of the quantum affine algebra
$U_q(\widehat{sl_N})$, in which the type $II$ VOs is not given
before. In fact, even for the classical affine Lie algebras, the
type $II$ VOs are not given. In this section, we will give the type
$II$ VOs of $U_q(\widehat{sl_N})$ and then use it to construct the
one of $U_{q,p}(\widehat{sl_N})$.

\vskip 0.5cm
\subsection{The type $I$ and type $II$ VOs of $U_q(\widehat{sl_N})$}
\vskip 0.5cm

In this subsection, we first review the primary field of the quantum
affine algebra $U_q(\widehat{sl_N})$ given in \cite{AOS}, which is
the type $I$ VOs of it. Here we denote it as
$\phi_{\vec{\Lambda}}(z)$ with $\vec{\Lambda}=(\lambda^1, \cdots,
\lambda^{N-1})$, where $\vec{\Lambda}$ is the weight of the
classical affine Lie algebra. However, we will reexpress it by using
some new bosons $\Big\{\check{a}^{i}: i=1, \cdots, N-1\Big\}$
defined as:
$$\check{a}^{i}_{n}=\sum_{j=1}^{N-1}n \frac{[\min(i, j)n]
[(N-\max(i, j))n]}{[(k+h^{\vee})n][Nn][n]^2}a_n^j$$
\no for any $n\in \mathbb{Z}_{\neq 0}$, and the zero-modes are

\beqas &&\check{p}_{a}^i=\sum_{j=1}^{N-1}\frac{\min(i, j)(N-\max(i,
j))}{(k+h^{\vee})N}p_a^j,\\
&&\check{q}_{a}^i=\sum_{j=1}^{N-1}\frac{\min(i, j)(N-\max(i,
j))}{(k+h^{\vee})N}q_a^j, \eeqas

\no here these bosons can be called the dual bosons of the original
ones $\Big\{a^{i}: i=1, \cdots, N-1\Big\}$ in the sense that they
satisfy the relations below: \beqas &&[\check{a}_n^i,
a_m^j]=\delta_{i,
j}\delta_{n+m, 0},\\
&&[\check{p}_a^i, q_{a}^j]=\delta_{i, j},\\
&&[\check{q}_a^i, p_{a}^j]=-\delta_{i, j},\eeqas

\no and these relations can be verified easily by using the q-analog
of the inverse of the Cartan matrix: \beqa
&&\sum_{r=1}^{N-1}\frac{[a_{i, r}n][\min(r, j)n][(N-\max(r,
j))n]}{[Nn][n]^2}=\delta_{i, j}.\eeqa

\no With them, the type $I$ VOs $\phi_{\vec{\Lambda}}(z)$ with
weight $\vec{\Lambda}=(\lambda^1, \cdots, \lambda^{N-1})$ of the
algebra $U_q(\widehat{sl_N})$ can be rewritten as:

$$\phi_{\vec{\Lambda}}(z)=:\prod_{i=1}^{N-1}\phi_{\lambda^{i}}(z):,$$

\no here the fields $\phi_{\lambda^{i}}(z)$ for $i=1,\cdots, N-1$,
which are called the components of $\phi_{\vec{\Lambda}}(z)$, are
given by:

\beqas
&&\phi_{\lambda^{i}}(z)=\exp\Bigg\{-\sum_{n>0}\frac{[\lambda^{i}n]}{n}
\check{a}^{i}_{-n}q^{\frac{k+h^{\vee}}{2}n}z^{n}\Bigg\}
\exp\Big\{\lambda^i(\check{q}_{a}^{i}+\check{p}_{a}^{i}\ln z)\Big\}
\exp\Bigg\{-\sum_{n>
0}\frac{[\lambda^{i}n]}{n}\check{a}_{n}^{i}q^{\frac{k+h^{\vee}}{2}n}z^{-n}\Bigg\};\eeqas

\no and it has the following properties with the fields
$\psi_{\pm}^{i}(z)$ and $e^{\pm, i}(z)$ in Eqns.(3.1)-(3.3) for
$i=1, \cdots, N-1$ proved in \cite{AOS} : \beqa
&&\psi_{\pm}^{i}(z)\phi_{\vec{\Lambda}}(w)=q^{-\lambda^{i}}\frac{w-q^{\lambda^{i}\mp
\frac{k}{2}}z}{w-q^{-\lambda^{i}\mp
\frac{k}{2}}z}\phi_{\vec{\Lambda}}(w)\psi_{\pm}^{i}(z);\\
&&e^{+, i}(z)\phi_{\vec{\Lambda}}(w)=\phi_{\vec{\Lambda}}(w)e^{+,
i}(z);\\
&&e^{-,
i}(z)\phi_{\vec{\Lambda}}(w)=q^{-\lambda^{i}}\frac{w-q^{\lambda^{i}
}z}{w-q^{-\lambda^{i}}z}\phi_{\vec{\Lambda}}(w)e^{-, i}(z).\eeqa

For any weight $\vec{\Lambda}=(\lambda^1, \cdots, \lambda^{N-1})$,
the type $II$ VOs $\psi_{\vec{\Lambda}}(z)$ of $U_q(\widehat{sl_N})$
is not known before. We present them in terms of $\Big\{b^{i, j},
c^{i, j}: 1\leq i< j\leq N-1\Big\}$ and the above mentioned dual
bosons. The field $\psi_{\vec{\Lambda}}(z)$ can be expressed as:
$$\psi_{\vec{\Lambda}}(z)=:\prod_{i=1}^{N-1}\psi_{\lambda^{i}}(z):$$

\no here its components $\psi_{\lambda^{i}}(z)$ are defined as:

\beqas &&\psi_{\lambda^{i}}(z)=\exp\Bigg\{-\sum_{n>0}
\frac{[\lambda^{N-i}n]}{n}\check{a}_{-n}^{i}q^{-\frac{k+h^{\vee}}{2}n}z^{n}\Bigg\}
\exp\Big\{\lambda^{N-i}(\check{q}_{a}^{i}+\check{p}_{a}^{i}\ln
z)\Big\}\exp\Bigg\{-\sum_{n>0}\frac{[\lambda^{N-i}n]}{n}
\check{a}_{n}^{i}q^{-\frac{k+h^{\vee}}{2}n}z^{-n}\Bigg\}\\
&&\hskip1.3cm \times \exp\Bigg\{\sum_{n>0}\sum_{j=i+1}^{N}
\frac{[\lambda^{j-i}n]}{[n]^2}(b_{-n}^{i,j}+c_{-n}^{i,
j})z^{n}\Bigg\}\exp\Bigg\{\sum_{j=i+1}^{N}\lambda^{j-i}\Big(q_{b}^{i,
j}+q_{c}^{i, j}+(p_{b}^{i, j}+p_{c}^{i,j})\ln
z\Big)\Bigg\}\\
&&\hskip1.3cm \times
\exp\Bigg\{-\sum_{n>0}\sum_{j=i+1}^{N}\frac{[\lambda^{j-i}n]}{[n]^2}(b_{n}^{i,
j}+c_{n}^{i, j})z^{-n}\Bigg\}.\eeqas

\no Furthermore, we have also proved the following theorem:

{\thm. The field $\psi_{\vec{\Lambda}}(z)$ satisfies the
intertwining relations: \beqas
&&\psi_{\pm}^{i}(z)\psi_{\vec{\Lambda}}(w)=q^{\lambda^{i}}
\frac{w-q^{-\lambda^{i}\pm
\frac{k}{2}}z}{w-q^{\lambda^{i}\pm
\frac{k}{2}}z}\psi_{\vec{\Lambda}}(w)\psi_{\pm}^{i}(z),\\
&&e^{+,
i}(z)\psi_{\vec{\Lambda}}(w)=q^{\lambda^{i}}
\frac{w-q^{-\lambda^{i}}z}{w-q^{\lambda^{i}}z}\psi_{\vec{\Lambda}}(w)e^{+,
i}(z),\\
&&e^{-, i}(z)\psi_{\vec{\Lambda}}(w)=\psi_{\vec{\Lambda}}(w)e^{-,
i}(z),\eeqas

\no where $\psi_{\pm}^{i}(z)$ and $e^{\pm, i}(z)$ for $ i=1, \cdots,
N-1$ are the currents given by (3.1)-(3.3).

\no Proof:}\ \ Here we only list the useful formulas we used to
prove this theorem: \beqas &&e^{A}B=e^{[A, B]}Be^{A},\ \ if\
[A, B]\ \ \ is\ \ a\ \ constant;\\
&&e^{A}e^{B}=e^{[A, B]}e^{B}e^{A}, \ \ \ if\ [A, B]\ \
commute\ \ with\ \ A\ \ and\ \ B;\\
&&\exp\Big(-\sum_{n>0}\frac{x^{n}}{n}\Big)=1-x;\\
&&(1-x)^{-1}=\sum_{n\geq 0}x^{n};\eeqas

\no and the q-analog of the inverse of the Cartan matrix in (5.22)
is also used. $\hskip0.5cm\Box$

These intertwining relations could be used to characterize the type
$II$ VOs of $U_q(\widehat{sl_N})$. Lastly, we present the
commutation relations among the type $I$ VOs
$\phi_{\vec{\Lambda}}(z)$ and type $II$ VOs
$\psi_{\vec{\Lambda}}(z)$; and here we only compute the commutation
relations between their components.

{\prop: \beqas
&&\phi_{\lambda^{i}}(z)\phi_{\lambda^{j}}(w)=\Big(\frac{z}{w}\Big)^{\lambda^{i}
\lambda^{j} g^{i,
j}}\exp\Big\{X_{1}\Big(\frac{z}{w}\Big)\Big\}\exp\Big\{-X_{1}(z\leftrightarrow
w)\Big\}
\phi_{\lambda^{j}}(w)\phi_{\lambda^{i}}(z);\\
&&\phi_{\lambda^{i}}(z)\psi_{\lambda^{j}}(w)=\Big(\frac{z}{w}\Big)^{\lambda^{i}
\lambda^{N-j} g^{i,
j}}\exp\Big\{X_{2}\Big(\frac{z}{w}\Big)\Big\}\exp\Big\{-X_{2}(z\leftrightarrow
w)\Big\}
\psi_{\lambda^{j}}(w)\phi_{\lambda^{i}}(z);\\
&&\psi_{\lambda^{i}}(z)\psi_{\lambda^{j}}(w)=\Big(\frac{z}{w}\Big)^{\lambda^{N-i}
\lambda^{N-j} g^{i,
j}}\exp\Big\{X_{3}\Big(\frac{z}{w}\Big)\Big\}\exp\Big\{-X_{3}(z\leftrightarrow
w)\Big\}\psi_{\lambda^{j}}(w)\psi_{\lambda^{i}}(z),\eeqas

\no here the $X_{i}\Big(\frac{z}{w}\Big)$ for $i=1, 2, 3$ are given
by

\beqas &&X_{1}\Big(\frac{z}{w}\Big)=\sum_{n>
0}\frac{1}{n}\frac{[\lambda^{i}n][\lambda^{j}n]}{[n]^2}[g^{i,
j}_{n}]q^{(k+h^{\vee})n}\Big(\frac{z}{w}\Big)^n;\\
&&X_{2}\Big(\frac{z}{w}\Big)=\sum_{n>
0}\frac{1}{n}\frac{[\lambda^{i}n][\lambda^{N-j}n]}{[n]^2}[g^{i,
j}_{n}]\Big(\frac{z}{w}\Big)^n;\\
&&X_{3}\Big(\frac{z}{w}\Big)=\sum_{n>
0}\frac{1}{n}\frac{[\lambda^{N-i}n][\lambda^{N-j}n]}{[n]^2}[g^{i,
j}_{n}]q^{-(k+h^{\vee})n}\Big(\frac{z}{w}\Big)^n,\eeqas

\no where for simplicity we use the symbols $g^{i, j}$ and $[g^{i,
j}_{n}]$ to denote:

\beqa &&g^{i, j}=\frac{\min(i, j)(N-\max(i, j))}{(k+h^{\vee})N},\\
&&[g^{i, j}_{n}]=\frac{[\min(i, j)n][(N-\max(i,
j))n]}{[(k+h^{\vee})n][N n]}. \eeqa }

\no It should be remarked that the matrix $G=\Big(g^{i,
j}\Big)_{1\leq i, j\leq N-1}$ is the inverse matrix of
$(k+h^{\vee})A$, and here $A$ is the Cartan matrix.

\vskip 0.5cm
\subsection{The type $I$ and type $II$ VOs of $U_{q,p}(\widehat{sl_N})$}
\vskip 0.5cm

In this subsection, we will give the free field realization of the
type $I$ and type $II$ VOs of $U_{q,p}(\widehat{sl_N})$. We nominate
them as $\Phi_{\vec{\Lambda}}(u)$ and $\Psi_{\vec{\Lambda}}(u)$ with
weight $\vec{\Lambda}=(\lambda^1, \cdots, \lambda^{N-1})$. They are
all obtained by twisting the corresponding ones of
$U_q(\widehat{sl_N})$ given in the above subsection.

First, we will construct two twisted currents $T_{\pm}(z; p)$ which
depend on the parameter $p$ for the two types VOs of
$U_q(\widehat{sl_N})$. For the type $I$ VOs, we define the twisted
current $T_{+}(z; p)$ as

$$T_{+}(z; p)=:\prod_{i=1}^{N-1}T_{+}^{i}(z; p):$$

\no and for $i=1, \cdots, N-1$,

$$T_{+}^{i}(z;
p)=\exp\Bigg\{\sum_{n>0}\frac{1}{n}\frac{[\lambda^{i}n][kn]}{[rn]}
\check{H}_{n}^{i}q^{(r-\frac{k}{2})n}z^{-n}\Bigg\}
\exp\Bigg\{-\frac{\lambda^{i}}{r}
\Big(\check{p}_{i}+\check{h}_{i}-(k+h^{\vee})\check{p}_{a}^{i}\Big)\ln
z\Bigg\}$$ \no here

\beqa &&\check{H}_{n}^{i}=\sum_{j=1}^{N-1}n \frac{[\min(i,
j)n][(N-\max(i, j))n]}{[kn][Nn][n]^2}H_{n}^{j},\ \ \ \ \forall n\in
\mathbb{Z}_{\neq 0}\eeqa

\no and

$$H_n^i=\sum_{j=1}^{i}\Big(b_{n}^{j,
i+1}q^{-(\frac{k}{2}+j-1)|n|}-b_{n}^{j,
i}q^{-(\frac{k}{2}+j)|n|}\Big)+a_{n}^{i}q^{-\frac{h^{\vee}}{2}|n|}
+\sum_{j=i+1}^{N}\Big(b_{n}^{i,
j}q^{-(\frac{k}{2}+j)|n|}-b_{n}^{i+1,
j}q^{-(\frac{k}{2}+j-1)|n|}\Big),$$ \no then it is easy to verify
that the following relation holds:
$$[\check{H}_{n}^{i}, H_m^j]=\delta_{i,j}\delta_{n+m, 0},$$

\no since we have the following commutation relation

$$[H_n^i, H_m^j]=\frac{[a_{ij}n][kn]}{n}\delta_{n+m, 0};$$

\no moreover, the symbols $\check{h}_{i}$ and $\check{p}_{i}$ are
used to denote the following complicated ones:

\beqas &&\check{h}_{i}=\sum_{j=1}^{N-1}\frac{\min(i,
j)(N-\max(i, j))}{N}h_{j},\\
&&\check{p}_{i}=\sum_{j=1}^{N-1}\frac{\min(i, j)(N-\max(i,
j))}{N}\hat{p}_{j}, \eeqas

\no and $\check{p}_{a}^i$ was defined at the beginning of the above
subsection. And for the type $II$ VOs, the twisted current $T_{-}(z;
p)$ is given by

$$T_{-}(z; p)=:\prod_{i=1}^{N-1}T_{-}^{i}(z; p):$$

\no and

$$T_{-}^{i}(z; p)=\exp\Bigg\{\sum_{n>0}\frac{1}{n}\frac{[\lambda^{i}n][kn]}
{[r^{*}n]}\check{H}_{-n}^{i}q^{(r-\frac{k}{2})n}z^{n}\Bigg\}
\exp\Big\{-2\lambda^{i}\check{q}_{i}
+\frac{\lambda^{i}}{r^{*}}\check{p}_{i}\ln z\Big\},$$

\no here $\check{H}_{-n}^{i}$ for $n>0$ is given in (5.28) and the
symbol $\check{q}_{i}$ is defined as follow:

$$\check{q}_{i}=\sum_{j=1}^{N-1}\frac{\min(i, j)(N-\max(i,
j))}{N}\hat{q}_{j}.$$

Next, we can use the components $T_{\pm}^{i}(z; p)$ of the twisted
currents to twist the ones $\phi_{\lambda^{i}}(z)$ and
$\psi_{\lambda^{i}}(z)$ as follows:

\beqas &&\Phi_{\lambda^{i}}(u)=:\phi_{\lambda^{i}}(z)T_{+}^{i}(z;
p):,\\
&&\Psi_{\lambda^{i}}(u)=:T_{-}^{i}(z; p)\psi_{\lambda^{i}}(z):,
\eeqas

\no here $\Phi_{\lambda^{i}}(u)$ and $\Psi_{\lambda^{i}}(u)$ could
also be considered as the components of the fields
$\Phi_{\vec{\Lambda}}(u)$ and $\Psi_{\vec{\Lambda}}(u)$, since we
define $\Phi_{\vec{\Lambda}}(u)$ and $\Psi_{\vec{\Lambda}}(u)$ by:

\beqas
&&\Phi_{\vec{\Lambda}}(u)=:\prod_{i=1}^{N-1}\Phi_{\lambda^{i}}(u): ,\\
&&\Psi_{\vec{\Lambda}}(u)=:\prod_{i=1}^{N-1}\Psi_{\lambda^{i}}(u):.
\eeqas

\no Furthermore, we obtain an important theorem by applying the
relations (5.23)-(5.25) and the Theorem $4$:

{\thm. The fields $\Phi_{\vec{\Lambda}}(u)$ and
$\Psi_{\vec{\Lambda}}(u)$ with given weight
$\vec{\Lambda}=(\lambda^1, \cdots, \lambda^{N-1})$ possess the
intertwining properties:

\beqa
&&H_i^{\pm}(u)\Phi_{\vec{\Lambda}}(v)=\frac{\theta_{r}(u-v
+\frac{\lambda^{i}}{2}\mp
\frac{k}{4})}
{\theta_{r}(u-v-\frac{\lambda^{i}}{2}\mp\frac{k}{4})}
\Phi_{\vec{\Lambda}}(v)H_i^{\pm}(u),\\
&&E_i(u)\Phi_{\vec{\Lambda}}(v)=\Phi_{\vec{\Lambda}}(v)E_i(u),\\
&&F_i(u)\Phi_{\vec{\Lambda}}(v)=\frac{\theta_{r}(u-v+\frac{\lambda^{i}}{2})}
{\theta_{r}(u-v-\frac{\lambda^{i}}{2})}\Phi_{\vec{\Lambda}}(v)F_i(u);\\
&&H_i^{\pm}(u)\Psi_{\vec{\Lambda}}(v)=\frac{\theta_{r^*}
(u-v-\frac{\lambda^{i}}{2}\pm
\frac{k}{4})}{\theta_{r^*}(u-v+\frac{\lambda^{i}}{2}\pm
\frac{k}{4})}\Psi_{\vec{\Lambda}}(v)H_i^{\pm}(u),\\
&&E_i(u)\Psi_{\vec{\Lambda}}(v)=\frac{\theta_{r^*}(u-v-\frac{\lambda^{i}}{2})}
{\theta_{r^*}(u-v+\frac{\lambda^{i}}{2})}\Psi_{\vec{\Lambda}}(v)E_i(u),\\
&&F_i(u)\Psi_{\vec{\Lambda}}(v)=\Psi_{\vec{\Lambda}}(v)F_i(u), \eeqa

\no where $H_i^{\pm}(u)$, $E_i(u)$ and $F_i(u)$ $(i=1, \cdots, N-1)$
are the total currents in (3.19)-(3.21).}

\no The above relations (5.29)-(5.34) could be used to define the
VOs of the elliptic quantum algebra $U_{q,p}(\widehat{sl_N})$. As a
result, we actually gave the free field realization of the type $I$
and type $II$ VOs of $U_{q,p}(\widehat{sl_N})$ with given level $k$.
Lastly, we also investigate the commutation relations among the VOs
$\Phi_{\vec{\Lambda}}(u)$ and $\Psi_{\vec{\Lambda}}(u)$:

{\prop. For the components $\Phi_{\lambda^{i}}(z)$ and
$\Psi_{\lambda^{i}}(z)$ of the type $I$ VOs and type $II$ VOs, we
have the following relations: \beqas
&&\Phi_{\lambda^{i}}(z)\Phi_{\lambda^{j}}(w)=\Big(\frac{z}{w}\Big)^{\lambda^{i}
\lambda^{j} g^{i,
j}}\exp\Big\{Y_1\Big(\frac{z}{w}\Big)\Big\}\exp\Big\{-Y_{1}(z\leftrightarrow w)\Big\}
\Phi_{\lambda^{j}}(w)\Phi_{\lambda^{i}}(z),\\
&&\Phi_{\lambda^{i}}(z)\Psi_{\lambda^{j}}(w)=\Big(\frac{z}{w}\Big)^{\lambda^{i}
\lambda^{N-j} g^{i, j}}\exp\Big\{-\frac{\lambda^{i}(k+g)}{r}\Big(
\lambda^{j}g^{i, j}+C^{i, j}\Big)\ln
z\Big\}\\
&&\hskip2.5cm
\times\exp\Big\{\Big(X_{2}+Y_{2}+Y_{3}+Y_{4}\Big)
\Big(\frac{z}{w}\Big)\Big\}\exp\Big\{-X_{2}(z\leftrightarrow
w)\Big\}\Psi_{\lambda^{j}}(w)\Phi_{\lambda^{i}}(z),\\
&&\Psi_{\lambda^{i}}(z)\Psi_{\lambda^{j}}(w)=\Big(\frac{z}{w}\Big)^{\lambda^{N-i}
\lambda^{N-j} g^{i,
j}}\exp\Big\{\Big(X_{3}+Y_{5}+Y_{6}\Big)\Big(\frac{z}{w}\Big)\Big\}\nn
\\&&\hskip2.5cm
\times \exp\Big\{-\Big(X_{3}+Y_{5}+Y_{6}\Big)(z\leftrightarrow w;
i\leftrightarrow j)\Big\}\Psi_{\lambda^{j}}(w)\Psi_{\lambda^{i}}(z),
\eeqas

\no here $g^{i, j}$ is given in (5.26); and the symbols $C^{i, j}$,
$\Big\{Y_{i}\Big(\frac{z}{w}\Big): i=1,\cdots, 6\Big\}$ are used to
simplify the complicated ones given below:

$$C^{i, j}=-\sum_{l=j}^{N-1}\lambda^{l+1-j}g^{i, l}+\sum_{l=j+1}^{N-1}\lambda^{l-j}g^{i, l}
-\Big(\sum_{l=j+1}^{N}\lambda^{l-j}\Big)g^{i,
j}+\Big(\sum_{l=j+1}^{N}\lambda^{l-j}\Big)g^{i, j-1};$$

\beqas &&Y_1\Big(\frac{z}{w}\Big)=-\sum_{n>
0}\frac{1}{n}\frac{[\lambda^{i}n][\lambda^{j}n][(r-k-h^{\vee})n]}{[rn][n]^2}[g^{i,
j}_{n}]q^{(r+k+h^{\vee})n}\Big(\frac{z}{w}\Big)^n,\\
&&Y_2\Big(\frac{z}{w}\Big)=-\sum_{n>
0}\frac{1}{n}\frac{[\lambda^{i}n][\lambda^{j}n][(k+h^{\vee})n]}{[rn][n]^2}[g^{i,
j}_{n}]q^{(r-k)n}\Big(\frac{z}{w}\Big)^n,\\
&&Y_3\Big(\frac{z}{w}\Big)=-\sum_{n>
0}\frac{1}{n}\frac{[\lambda^{i}n][\lambda^{N-j}n][(k+h^{\vee})n]}{[rn][n]^2}[g^{i,
j}_{n}]q^{(r-k-h^{\vee})n}\Big(\frac{z}{w}\Big)^n,\\
&&Y_4\Big(\frac{z}{w}\Big)=-\sum_{n>
0}\frac{1}{n}\frac{[\lambda^{i}n][(k+h^{\vee})n]}{[rn][n]^2}[g^{i,
j}_{n}][C^{i, j}]q^{(r-k/2)n}\Big(\frac{z}{w}\Big)^n,\\
&&Y_5\Big(\frac{z}{w}\Big)=-\sum_{n>
0}\frac{1}{n}\frac{[\lambda^{N-i}n][\lambda^{j}n][(k+h^{\vee})n]}{[r^{*}n][n]^2}[g^{i,
j}_{n}]q^{(r-k-h^{\vee})n}\Big(\frac{z}{w}\Big)^n,\\
&&Y_6\Big(\frac{z}{w}\Big)=-\sum_{n>
0}\frac{1}{n}\frac{[\lambda^{j}n][(k+h^{\vee})n]}{[r^{*}n][n]^2}[g^{i,
j}_{n}]\Big([C^{i, j}](i\leftrightarrow
j)\Big)q^{(r-k/2)n}\Big(\frac{z}{w}\Big)^n ,\eeqas

\no where $[g^{i, j}_{n}]$ is given by (5.27) and $[C^{i, j}]$ is
defined by:

\beqas &&[C^{i, j}]=-\Big(\sum_{l=j}^{N-1}[\lambda^{l+1-j}n][g^{i,
l}_{n}]\Big)q^{-(k/2+j-1)n}+\Big(\sum_{l=j+1}^{N-1}[\lambda^{l-j}n][g^{i,
l}_{n}]\Big)q^{-(k/2+j)n}\\
&&\hskip 1.5cm
-\Big(\sum_{l=j+1}^{N}[\lambda^{l-j}n]q^{-(k/2+l)n}\Big)[g^{i,
j}_{n}]+\Big(\sum_{l=j+1}^{N}[\lambda^{l-j}n]q^{-(k/2+l-1)n}\Big)[g^{i,
j-1}_{n}].\eeqas}

\vskip 1cm
\section{Discussion}
\vskip 0.5cm

In this paper, we construct the free field representation of
$U_{q,p}(\widehat{sl_N})$ with given level $k$ by twisting the
Wakimoto realization of the quantum affine algebra
$U_q(\widehat{sl_N})$. The free boson realization of its screening
currents are also given. Moreover, the explicit expressions of the
type $II$ VOs of $U_q(\widehat{sl_N})$ and the two types VOs of
$U_{q,p}(\widehat{sl_N})$ are presented. In fact, even for the
classical affine Lie algebras, the type $II$ VOs are not given. We
also have much interests in the derivation of the multi-point
correlation functions, but in view of its complexity and the length
of the manuscript, it will be discussed in the future. Meanwhile, it
is also very interesting to extend our results to other types of Lie
algebras, and we will discuss them in a separate paper.

\vskip 1cm
\section{Acknowledgments}
\vskip 0.5cm

The results in section $3$ were reported in the XXth Conference of
Lie Algebras in Changshu, China (2007). One of the authors (Ding) is
financially supported partly by the Natural Science Foundations of
China through the grands No.10671196 and No.10231050. He is also
supported partly by a fund of Chinese Academy of Sciences.

\bebb{99}

\bbit{CFT} Di Francesco, P., Mathieu, P., S\'{e}n\'{e}chal, D.: {\it
Conformal Field Theory}. New York: Springer-Verlag, 1997

\bbit{Drinfeld86} Drinfeld, V. G.: Quantum groups. In {\it Proc.of
ICM Berkeley 1986}, Providence, RI: AMS, 798-820 (1987)

\bbit{Jimbo} Jimbo, M.: A q-difference analogue of $U~g$ and the
Yang-Baxter equation. Lett. Math. Phys. {\bf 10}, 62-69 (1986)

\bbit{Smir} Smirnov, F.A.: Dynamical symmetries of massive
integrable models $I$. Int. J. Mod. Phys. A {\bf 7} suppl.{\bf 1B},
813-837 (1992); Dynamical symmetries of massive integrable models
$II$. Int. J. Mod. Phys. A{\bf 7} suppl.{\bf 1B}, 839-858 (1992)

\bbit{Felder} Felder, G.: Elliptic quantum groups. In {\it Proc.
ICMP Pairs 1994}, Cambridge-Hong Kong: International Press, 211-218
(1995)

\bbit{Fronsdal} Fr{\o}nsdal, C.: Generalization and exact
deformation of quantum groups. Kyoto Univ.: Publ. RIMS, 91-149
(1997)

\bbit{EF} Enriquez, B., Felder, G.: Elliptic quantum groups
$E_{\tau, \eta}(sl_2)$ and quasi-Hopf algebras.  Commun. Math. Phys.
{\bf 195}, 651-689 (1998)

\bbit{Drinfeld} Drinfeld, V.G.: Quasi-Hopf algebras. Leningrad Math.
Jour.{\bf 1}, 1419-1457 (1990)

\bibitem{Kac} Kac, V.G.: {\it Infinite dimensional Lie algebras}, third ed..
 Cambridge: Cambridge University Press, 1990

\bbit{Baxter} Baxter, R.J.: Partition function of the eight-vertex
lattice model. Ann. Physics {\bf 70}, 193-228 (1972)

\bbit{ABF} Andrew, G.E., Baxter, R.J., Forrester, P.J.: Eight-vertex
SOS model and generalized Rogers-Ramanujan-type identities. J. Stat.
Phys. {\bf 35}, 193-266 (1984)

\bbit{FaddTak} Faddeev, L., Takhtajan, L.: {\it Hamiltonian methods
in the theory of solitons}. Berlin: Springer-Verlag, 1987

\bbit{JimMiw} Jimbo, M., Miwa, T.: {\it Algebraic analysis of
solvable lattice models}. CBMS Regional Conference Series in
Mathematics, Vol. {\bf 85}, Providence, RI: AMS, 1994

\bbit{KZ} Knizhnik, V.G., Zamolodchikov, A.B.: Current algebra and
Wess-Zumino models in two dimensions. Nucl. Phys. B{\bf 247}, 83-103
(1984)

\bbit{FR} Frenkel, I.B., Reshetikhin, N.Y.: Quantum affine algebras
and holomorphic difference equations. Commun. Math. Phys. {\bf 146},
 1-60 (1992)

\bbit{Bern} Bernard, D.: On the Wess-Zumino-Witten model on the
torus. Nucl. Phys. B{\bf 303}, 77-93 (1988)

\bibitem{Wak}%{}
 Wakimoto, M.: Fock representations of the affine Lie algebra
$A_{1}^{(1)}$.  Commun. Math. Phys. {\bf 104}, 605-609 (1986)

\bbit{FF} Feigin, B., Frenkel, E.: Affine Kac-Moody algebras and
semi-infinite flag manifolds.  Commun. Math. Phys. {\bf 128},
161-189 (1990)

\bbit{Yu} Petersen, J.L., Rasmussen, J., Yu, M.: Free field
realizations of 2D current algebras, screening currents and primary
fields. Nucl. Phys. B{\bf 502}, 649-670 (1997)

\bbit{Nemesch} Nemeschansky, D.:  Feigin-Fuchs representation of
$\widehat{su(2)}_k$ Kac-Moody algebra. Phys. Lett. B{\bf 224},
 121-124 (1989)

\bbit{GepQiu} Gepner, D., Qiu, Z.: Modular invariant partition
functions for parafermionic field theories. Nucl. Phys. B{\bf 285},
423-453 (1987)

\bbit{Gep} Gepner, D.: New conformal field theories associated with
Lie algebras and their partition functions. Nucl. Phys. B{\bf 290},
 10-24 (1987)

\bbit{Feh1}Boer, J., Feher, L.: Wakimoto realizations of current
algebras: an explicit construction. Commun.Math.Phys. \textbf{189},
759-793 (1997)

\bbit{Feh2}Feher, L., Pusztai, B.G.: Explicit description of twisted
Wakimoto realizations of affine Lie algebras. Nucl.Phys.
B\textbf{674}, 509-532 (2003)

\bbit{YangZh}Yang, W.-L., Zhang, Y.-Z.: On explicit free fields
realizations of current algebras. Nucl.Phys.B.{\bf 800}, 527-546
(2008)

\bbit{DFJMN} Davies, B., Foda, O., Jimbo, M., Miwa, T., Nakayashiki,
A.: Diagonalization of the XXZ Hamiltonian by vertex operators.
Commun.Math.Phys. {\bf 151}, 89-153 (1993)

\bbit{JMMN} Jimbo, M., Miki, K., Miwa, T., Nakayashiki, A.:
Correlation functions of the XXZ model for $\triangle <{-1}$.  Phys.
Lett. A{\bf 168}, 256-263 (1992)

\bbit{Shiraishi} Shiraishi, J.: Free boson representation of
$U_q(\widehat{sl_2})$. Phys, Lett. A{\bf 171}, 243-248 (1992)

\bbit{Matsuo} Matsuo, A.: Free field representation of quantum
affine algebra $U_q(\widehat{sl_2})$.  Phys. Lett. B{\bf 308},
 260-265 (1993)

\bbit{DingW} Ding, X.-M., Wang, P.: Parafermion representation of
the quantum affine $U_q(\widehat{sl_2})$. Modern. Phys. Lett. A{\bf
11}, 921-930 (1996)

\bbit{AOS} Awata, H., Odake, S., Shiraishi, J.: Free boson
realization of $U_q(\widehat{sl_N})$.  Commun.Math.Phys. {\bf 162},
61-83 (1994)

\bbit{Luky} Lukyanov, S.: Free field representation for massive
integrable models. Commun.Math.Phys. {\bf 167}, 183-226 (1995)

\bbit{IK} Iohara, K., Kohno, M.: A central extension of
$DY_{\hbar}(\widehat{sl_2})$ and its vertex representation.  Lett.
Math. Phys. {\bf 37}, 319-328 (1996)

\bbit{Konno2} Konno, H.: Free field representation of level-$k$
Yangian double $DY_{\hbar}(\widehat{sl_2})_{k}$ and deformation of
Wakimoto modules. Lett. Math. Phys. {\bf 40}, 321-336 (1997)

\bbit{Iohara} Iohara, K.: Bosonic representations of Yangian double
$DY_{\hbar}(g)$ with $g=gl_N$, $sl_N$.  J. Phys. A: Math. Gen.{\bf
29}, 4593-4621 (1996)

\bbit{DHHZ}Ding, X.-M., Hou, B.-Y., Hou, B.-Yuan, Zhao, L.: Free
boson representation of $DY_{\hbar}({gl_N})_{k}$ and
$DY_{\hbar}({sl_N})_{k}$.  J. Math. Phys. {\bf 39}, 2273-2289 (1998)

\bbit{DHZ}Ding, X.-M., Hou, B.-Y., Zhao, L.: $\hbar$-Yangian
deformation of the Miura map and Virasoro algebra.  Intern. Jour.
Mod. Phys. A{\bf 13}, 1129-1144 (1998)

\bbit{HY} Hou, B.-Y., Yang, W.-L.: A $\hbar$-deformed Virasoro
algebra as hidden symmetry of the restricted sine-Gordon model.
Commun. Theor. Phys. {\bf 31}, 265-270 (1999)

\bbit{Lukpugai} Lukyanov, S., Pugai, Y.: Multi-point Local Height
Probabilities in the Integrable RSOS Model. Nucl.Phys. B{\bf 473},
631-658 (1996)

\bbit{Konno1} Konno, H.: An elliptic algebra
$U_{q,p}(\widehat{sl_2})$ and the fusion RSOS model.  Commun. Math.
Phys. {\bf 195}, 373-403 (1998)

\bbit{CD07} Chang, W.-J., Ding, X.-M.: On the vertex operators of
the elliptic quantum algebra $U_{q,p}(\widehat{sl_2})_{k}$.  J.
Math. Phys. {\bf 49}, 043513 (2008)

\bbit{JKOS} Jimbo, M., Konno, H., Odake, S., Shiraishi, J.: Elliptic
algebra $U_{q,p}(\widehat{sl_2})$: Drinfeld currents and vertex
operators. Commun. Math. Phys. {\bf 199}, 605-647 (1999)

\bbit{KK} Kojima, T., Konno, H.: The elliptic algebra
$U_{q,p}(\widehat{sl_N})$ and the Drinfeld realization of the
elliptic quantum group $B_{q,\lambda}(\widehat{sl_N})$. Commun.
Math. Phys. {\bf 239}, 405-447 (2003)

\bbit{Dong} Dong, C.-Y., Lepowsky, J.: {\it Generalized Vertex
Algebras and Relative Vertex Operators}, Prog. Math.{\bf 112}. Boston:
Birh\"{a}user, 1993

\eebb

\end{document}